# COMPLETE CORRECTED DIFFUSION APPROXIMATIONS FOR THE MAXIMUM OF A RANDOM WALK

BY JOSE BLANCHET AND PETER GLYNN

*Harvard University and Stanford University*

Consider a random walk $(S_n : n \geq 0)$ with drift $-\mu$ and $S_0 = 0$. Assuming that the increments have exponential moments, negative mean, and are strongly nonlattice, we provide a complete asymptotic expansion (in powers of $\mu > 0$) that corrects the diffusion approximation of the all time maximum $M = \max_{n \geq 0} S_n$. Our results extend both the first-order correction of Siegmund [*Adv. in Appl. Probab.* **11** (1979) 701–719] and the full asymptotic expansion provided in the Gaussian case by Chang and Peres [*Ann. Probab.* **25** (1997) 787–802]. We also show that the Cramér–Lundberg constant (as a function of $\mu$) admits an analytic extension throughout a neighborhood of the origin in the complex plane $\mathbb{C}$. Finally, when the increments of the random walk have nonnegative mean $\mu$, we show that the Laplace transform, $E_\mu \exp(-bR(\infty))$, of the limiting overshoot, $R(\infty)$, can be analytically extended throughout a disc centered at the origin in $\mathbb{C} \times \mathbb{C}$ (jointly for both $b$ and $\mu$). In addition, when the distribution of the increments is continuous and appropriately symmetric, we show that $E_\mu S_\tau$ [where $\tau$ is the first (strict) ascending ladder epoch] can be analytically extended to a disc centered at the origin in $\mathbb{C}$, generalizing the main result in [*Ann. Probab.* **25** (1997) 787–802] and extending a related result of Chang [*Ann. Appl. Probab.* **2** (1992) 714–738].

**1. Introduction.** Let $(X_n : n \geq 1)$ be a sequence of independent and identically distributed (i.i.d.) random variables (r.v.'s), and let $S = (S_n : n \geq 0)$ be its associated random walk (so that $S_0 = 0$ and $S_n = X_1 + \cdots + X_n$ for $n \geq 1$). This paper is a contribution to the literature on development of high accuracy approximations to the distribution of the maximum r.v.

$$M = \max\{S_n : n \geq 0\}.$$

Clearly, $-\mu \triangleq EX_1$ must typically be negative in order that $M$ be finite-valued. We shall require this condition throughout the remainder of this paper. The distribution of $M$ is of importance in a number of different disciplines.

For $x > 0$, $\{M > x\} = \{\tau(x) < \infty\}$, where $\tau(x) = \inf\{n \geq 1 : S_n > x\}$, so that computing the tail of $M$ is equivalent to computing a level crossing probability for the random walk $S$. Because of this level crossing interpretation, the tail of $M$ is of great interest to both the sequential analysis and risk theory communities. In particular, in the setting of insurance risk, $P(\tau(x) < \infty)$ is the probability that an insurer will face ruin in finite time (when the insurer starts with initial reserve $x$ and is subjected to i.i.d. claims over time); see, for example, [1].

The distribution of $M$ also arises in the analysis of the single most important model in queueing theory, namely, the single-server queue. If the inter-arrival and service times for successive customers are i.i.d. with a mean arrival rate less than the mean service rate, then $W = (W_n : n \geq 0)$ is a positive recurrent Markov chain on $[0, \infty)$, where $W_n$ is the waiting time (exclusive of service) for customer $n$. If $W_\infty$ is a random variable having the stationary distribution of $W$, then Kiefer and Wolfowitz [9] showed that $W_\infty$ has the distribution of $M$ for an appropriately defined random walk. As a consequence, computing the distribution of $M$ is of fundamental importance to queueing theorists.

Since $W$ is a positive recurrent Markov chain, the distribution of $M$ can be computed as the solution to the equation describing the stationary distribution of $W$. This linear integral equation is known as Lindley's equation (see [12]) and is of Wiener–Hopf type; it is challenging to solve, both analytically and numerically. As a result, approximations are frequently employed instead. One important such approximation holds as $\mu \searrow 0$. This asymptotic regime corresponds in risk theory to the setting in which the "safety loading" is small (i.e., the premium charged is close to the typical pay-out for claims) and in queueing theory to the "heavy traffic" setting in which the server is utilized close to 100% of the time. Thus, this asymptotic regime is of great interest from an applications standpoint. Kingman [10] showed that the approximation

$$P(M > x) \approx \exp(-2\mu x/\sigma^2) \tag{1}$$

is valid as $\mu \searrow 0$, where $\sigma^2 = \text{Var}(X_1)$. (A more precise statement of this result will be given in Section 2.) Because the right-hand side of (1) is the exact value of the level crossing probability for the natural Brownian approximation to the random walk $S$, (1) is often called the diffusion approximation to the distribution of $M$.

As with any such approximation, there are applications for which (1) delivers poor results. Siegmund [15] therefore proposed a so-called "corrected diffusion approximation" that reflects information in the increment distribution beyond the mean and variance. This corrected diffusion approximation computes the next term in the asymptotic (as $\mu \searrow 0$) beyond that given by the right-hand side of (1). The main result in this paper (Theorem 1) is a development of the full asymptotic expansion initiated by Siegmund. We compute all the terms in the asymptotic expansion for general random walks with increments having exponential moments; see Section 6 for details on the calculation of the relevant coefficients in the expansion. Our theorem can be viewed as a non-Gaussian counterpart to the corresponding expansion provided recently by Chang and Peres [7] for Gaussian random walks. As perhaps expected, the mathematical approach followed here is quite different from that used by Chang and Peres.

In addition to risk theory and queueing theory, the types of results and techniques studied here can be useful in other settings such as statistical sequential analysis, which is the discipline that motivated the initial development of Siegmund [15]; see also [16]. Other areas that have benefited from corrected diffusion approximations also include inventory theory; see [8] and computational finance; see, for example, [3].

As is well known in the literature, there is a close connection between such corrections and asymptotic expansions for the moments of the ascending ladder height random variables associated with the random walk. Theorem 2 establishes an asymptotic expansion for the mean of the first strict ascending ladder height for random walks with light-tailed symmetric and continuous increments. As indicated in Section 6, this permits one to develop asymptotic expansions for all the moments of the ascending ladder heights (and for the limiting overshoot induced by the associated renewal process); see also Theorem 4. Additional expressions on boundary crossing problems and ladder height quantities are obtained by Lai [11] and, in a different setting in which $\mu$ is held fixed, by Lotov [13].

This paper is organized as follows. The main results are described in Section 2. A key connection to asymptotic expansions for the "short-time" behavior of the Cauchy process is made in Section 3. Section 4 shows how all the integrals required for our asymptotic expansion can be reduced to the short-time asymptotics of Section 3. Finally, Section 5 provides rigorous support for the remaining details in the argument used to compute the coefficients in the expansion. Section 6 summarizes the computation of the coefficients, and discusses an expansion related to the moments of the strict ascending ladder height. Any proof that does not follow the statement of the result can be found in our final section, namely, Section 7.

**2. The main results.** To state our main results, we adopt the parameterization utilized by Siegmund [15]. We assume throughout this paper that the $X_i$'s have exponential moments, so that $E\exp(\theta X_1) < \infty$ for $\theta$ in a neighborhood containing the origin. For such $\theta$, define

$$\psi(\theta) = \log E(\exp(\theta X_1)).$$

Then, for each such $\theta$, we can define the probability measure $P_\theta$ having the property that, for $n \geq 0$,

$$P_\theta(A) = E(\exp(\theta S_n - n\psi(\theta))\mathbb{1}_A)$$

for $A \in \sigma(S_j : 0 \leq j \leq n)$. As is well known, $S$ is again a random walk with i.i.d. increments under $P_\theta$, having common increment distribution

$$P_\theta(X_1 \in dx) = \exp(\theta x - \psi(\theta))P(X_1 \in dx)$$

for $x \in \mathbb{R}$ [with mean $E_\theta X_1 = \psi'(\theta)$ and variance $\mathrm{Var}_\theta(X_1) = \psi''(\theta)$]. Without any loss of generality, assume that $EX_1 = 0$ and $\mathrm{Var}(X_1) = 1$. Since $\psi(\cdot)$ is strictly convex on its domain of finiteness, $E_\theta X_1 < 0$ for $\theta < 0$. Thus, $P_\theta$ induces a random walk with negative drift when $\theta < 0$. We therefore focus on corrected approximations to $P_\theta(M > x)$ as $\theta \nearrow 0$.

A key step to the analysis of $P_\theta(M > x)$ is the judicious application of Wald's likelihood ratio identity; see, for example, [16], page 13. For $\theta_0$ in some interval of the form $(-\eta, 0)$, there exists a positive $\theta_1$ such that $\psi(\theta_0) = \psi(\theta_1)$. Set $\Delta = \theta_1 - \theta_0$. Note that parameterizing in terms of $\Delta$ is essentially equivalent to parameterization in terms of $\theta_0$ [or parameterization in terms of the drift $\mu = -\psi'(\theta_0)$]. The likelihood ratio identity then asserts that

$$\begin{aligned}
(2)\quad P_{\theta_0}(\tau(x) < \infty) &= E_{\theta_1}\exp(-(\theta_1 - \theta_0)S_{\tau(x)}) \\
&= \exp(-(\theta_1 - \theta_0)x)E_{\theta_1}\exp(-(\theta_1 - \theta_0)R(x)),
\end{aligned}$$

where $R(x) = S_{\tau(x)} - x$ is the so-called "overshoot" at level $x$.

Suppose now that $X_1$ is strongly nonlattice, in the sense that, for each $\delta > 0$,

$$(3) \quad \inf_{|\lambda| > \delta} |1 - g(\lambda)| > 0,$$

where $g(\lambda) = E\exp(i\lambda X_1)$ is the characteristic function of $X_1$ (under $P_0$). Applying renewal theory to the random walk at strictly increasing ladder epochs establishes then

$$(4) \quad E_{\theta_1}\exp(-(\theta_1 - \theta_0)R(x)) \to E_{\theta_1}\exp(-(\theta_1 - \theta_0)R(\infty))$$

as $x \to \infty$.

Siegmund [15] showed that the renewal theorem can be applied uniformly for $\Delta < \eta$ (see also [6]). Hence, (2) yields

$$(5) \quad P_{\theta_0}(M > x) = \exp(-\Delta x)E_{\theta_1}\exp(-\Delta R(\infty)) + o(\exp(-(\Delta + r)x))$$

for some $r > 0$ (uniformly in $\theta_0 > -\eta/2$). In insurance risk theory, $\Delta$ is called the "adjustment coefficient" and the quantity $E_{\theta_1} \exp(-\Delta R(\infty))$ is known as the Cramér–Lundberg constant [1].

Relation (5) may alternatively be written as

(6) $\quad P_{\theta_0}(\Delta M > x) = \exp(-x) E_{\theta_1} \exp(-\Delta R(\infty)) + o(\exp(-rx/\Delta))$,

where $o(\exp(-rx/\Delta))$ is uniform in $\theta_0 > -\eta/2$. Note that $\exp(-x)$ is precisely the level crossing probability of level $x/\Delta$ for a Brownian motion with drift $-\Delta/2$ and unit variance. Since $E_{\theta_1} X_1 \sim -\Delta/2$ as $\theta_0 \nearrow 0$, (6) provides rigorous support for the diffusion approximation (1). Furthermore, a correction to the diffusion approximation described in the Introduction can be obtained by developing an asymptotic expansion for $E_{\theta_1} \exp(-\Delta R(\infty))$.

Siegmund [15] obtained his corrected diffusion approximation by showing that

(7) $\quad\quad\quad\quad E_{\theta_1} \exp(-\Delta R(\infty)) = \exp(-\Delta \beta_1) + o(\Delta^2)$

as $\Delta \downarrow 0$, where $\beta_1$ can be computed explicitly as

(8) $\quad\quad\quad \beta_1 = \frac{1}{6} E X_1^3 - \frac{1}{2\pi} \int_{-\infty}^{\infty} \frac{1}{\lambda^2} \operatorname{Re} \log\{2(1-g(\lambda))/\lambda^2\} \, d\lambda.$

Note that, by computing the single integral (8), Siegmund's corrected diffusion approximation to the distribution of $M$ provides a parametric approximation that is valid for all random walks having negative drift sufficiently close to zero. Such parametric approximations are convenient in many applications settings (i.e., in studying the behavior of a queue when utilization is close to 100%).

Our main theorem shows that there is a full asymptotic expansion for

$$r(\Delta) \triangleq \log E_{\theta_1} \exp(-\Delta R(\infty)).$$

THEOREM 1. *Suppose that $X_1$ has exponential moments and is strongly non-lattice. Then, $r(\cdot)$ (initially defined on $[0, \upsilon)$ for some $\upsilon > 0$) admits an analytic extension on a neighborhood of the origin in the complex plane.*

REMARK. An immediate consequence of Theorem 1 and the implicit function theorem is that the Cramér–Lundberg constant, namely, $\exp(r(\Delta(\theta_0)))$, initially defined for all $\theta_0 < 0$ sufficiently close to zero, admits an analytic extension on a disc containing the origin in the complex plane.

According to Theorem 1,

(9) $\quad\quad\quad\quad E_{\theta_1} \exp(-\Delta R(\infty)) = \exp\left(\sum_{n=1}^{\infty} \beta_n \Delta^n\right),$

where $\beta_1$ is given by (8) and $\beta_2 = 0$. [This latter equality follows from the fact that the error term in (7) is $o(\Delta^2)$.] Obviously, in order for (9) to be useful from an applied standpoint, we need a means of numerically computing the $\beta_n$'s. This issue is discussed in Section 6. We establish there that the $\beta_n$'s can be successively computed via a finite number of one-dimensional integrations reminiscent of the integral appearing in (8). Thus, the $\beta_n$'s can easily be computed, thereby yielding cheaply computable high-order parametric corrections to the diffusion approximation (1).

The argument above also permits us to establish asymptotic expansions for certain ladder height quantities. As noted earlier, renewal theory applies to the random walk when sampled at strictly increasing ladder epochs. The renewal theorem invoked above actually establishes that

$$(10) \qquad E_{\theta_1} \exp(-\Delta R(\infty)) = \frac{1 - E_{\theta_1} \exp(-\Delta S_{\tau_+})}{\Delta E_{\theta_1} S_{\tau_+}},$$

where $\tau_+ = \inf\{n \geq 1 : S_n > 0\}$ is the first (strict) increasing ladder epoch (see [2], page 221). In view of (2), it follows that

$$(11) \qquad 1 - E_{\theta_1} \exp(-\Delta S_{\tau_+}) = P_{\theta_0}(\tau_+ = \infty).$$

Random walk duality (see, e.g., page 173 of [16]) implies that

$$(12) \qquad P_{\theta_0}(\tau_+ = \infty) = 1/E_{\theta_0}\tau_-,$$

where $\tau_- = \inf\{n \geq 1 : S_n \leq 0\}$. If the $X_i$'s are symmetric r.v.'s with common continuous distribution function, $\Delta = 2\theta_1$ and $E_{\theta_0}\tau_- = E_{\theta_1}\tau_+$. Furthermore, (10) to (12) then imply that

$$E_{\theta_1} \exp(-\Delta R(\infty)) = \frac{1}{2\theta_1 (E_{\theta_1} S_{\tau_+})(E_{\theta_1}\tau_+)}.$$

In view of Wald's identity, we then obtain the relation

$$E_{\theta_1} \exp(-\Delta R(\infty)) = \frac{\psi'(\theta_1)}{2\theta_1 (E_{\theta_1} S_{\tau_+})^2}.$$

As a consequence, Theorem 1 then yields a full asymptotic expansion for the expected ladder height $E_{\theta_1} S_{\tau_+}$. We record this result as our Theorem 2.

THEOREM 2. *Assume that $X_1$ has exponential moments and is symmetric with a continuous distribution function. Then,*

$$E_{\theta_1} S_{\tau_+} = \sqrt{\frac{\psi'(\theta_1)}{2\theta_1}} \exp\left(-\frac{1}{2} \sum_{m=0}^{\infty} \beta_{2m+1}(2\theta_1)^{2m+1}\right).$$

Given our above argument, the only remaining issue in proving Theorem 2 is establishing that $\beta_{2n} = 0$ for $n \geq 1$ in the presence of symmetry. This fact is proven in Section 7.

The most important device that we use to prove Theorems 1 and 2 is a convenient representation for $r(\Delta)$. This representation, is a key idea in our mathematical development. To introduce our representation, put $\phi(\theta) = E\exp(\theta X_1)$ for $\theta \in \mathbb{R}$ and, for $z \in \mathbb{C}$, set $\gamma(z) = E\exp(zX_1)$. Note that $\phi$ is finite-valued on a neighborhood $\mathcal{N}$ of the origin and $\gamma$ is analytic on the strip $\{x + iy : x \in \mathcal{N}, y \in \mathbb{R}\}$. For nonnegative $\theta \in \mathcal{N}$ and $b \in \mathbb{R}$, put

$$\rho(\theta, b) = \log E_\theta \exp(-bR(\infty)).$$

Note $r(\Delta) = \rho(\theta_1, \Delta)$, where $\theta_1 = \theta_1(\Delta) > \theta_0(\Delta) = \theta_0$ is such that $\psi(\theta_1(\Delta)) = \psi(\theta_0(\Delta))$. Woodroofe [17] showed that

$$(13) \qquad \rho(\theta, b) = \frac{1}{2\pi} \int_{-\infty}^{\infty} \frac{-b}{(b+i\lambda)i\lambda} \log\left(\frac{\gamma(\theta) - \gamma(\theta + i\lambda)}{-i\phi'(\theta)\lambda}\right) d\lambda;$$

see also Corollary 8.45 and Theorem 8.51 of [16]. While (13) is convenient for many purposes, it presents difficulties in the current circumstances because of the singularity (in the logarithm) that arises when $\theta \searrow 0$. The following representation for $\rho(\theta, b)$ is free of such singularities.

THEOREM 3. *Suppose $X_1$ has exponential moments and is strongly nonlattice. Then, for nonnegative $\theta \in \mathcal{N}$ and $b > 0$,*

$$(14) \qquad \rho(\theta, b) = \frac{1}{2\pi} \int_{-\infty}^{\infty} \frac{-b}{(b+i\lambda)i\lambda} \log\left(\frac{2(\gamma(\theta) - \gamma(\theta + i\lambda))}{\lambda(\lambda - 2i\phi'(\theta))}\right) d\lambda.$$

Siegmund's computation of $\beta_1$ takes advantage of the fact that the first-order behavior of $r(\Delta)$ should match that of

$$(15) \qquad s(\Delta) = \log E_0 \exp(-\Delta R(\infty)).$$

Since $s(\Delta) = \rho(0, \Delta)$, Theorem 3 implies that

$$(16) \qquad s(\Delta) = \frac{1}{2\pi} \int_{-\infty}^{\infty} \frac{-\Delta}{(\Delta + i\lambda)i\lambda} \log(2(1 - g(\lambda))\lambda^{-2}) d\lambda;$$

see also page 226 of [16]. We proceed to analyze $\rho(\theta, b)$ by writing $\rho(\theta, b) = s(b) + I(\theta, b)$. In view of both Theorem 3 and (16),

$$(17) \qquad I(\theta, b) = \frac{1}{2\pi} \int_{-\infty}^{\infty} \frac{-b}{(b+i\lambda)i\lambda} \log\left(\frac{\lambda(\gamma(\theta) - \gamma(\theta + i\lambda))}{(\lambda - 2i\phi'(\theta))(1 - g(\lambda))}\right) d\lambda.$$

In the next sections we develop asymptotics, as $b \searrow 0$, appropriate to the integrals arising in (16) and (17). Such asymptotics can be used to provide asymptotic expansions for the moments (or, equivalently, the cumulants) of the limiting expected overshoot r.v. $R(\infty)$ under $P_\theta$ as $\theta \searrow 0$. Specifically, for $n \geq 1$, let

$$\kappa_n(\theta) = (-1)^n \frac{\partial^n}{\partial b^n} \rho(\theta, b) \bigg|_{b=0}.$$

THEOREM 4. *Assume that $X_1$ has exponential moments and is strongly non-lattice. Then (for all $n \geq 1$) $\kappa_n(\cdot)$, initially defined on $[0, \upsilon)$ for $\upsilon > 0$, can be extended to be an analytic function throughout a disc in the complex plane containing the origin.*

An important implication of Theorem 4 is that it can be directly applied to obtain complete asymptotics for the steady-state mean of the waiting time sequence, namely, $E_{\theta_0} M \, (= E_{\theta_0} W_\infty)$. In particular, Siegmund [15] shows (see also page 270 of [2]) that

(18)
$$\begin{aligned}
E_{\theta_0} M &= \frac{E_{\theta_0}(S_{\tau_+} | \tau_+ < \infty)}{P_{\theta_0}(\tau_+ = \infty)} \\
&= \frac{E_{\theta_1} S_{\tau_+} \exp(-\Delta S_{\tau_+})}{1 - E_{\theta_1} \exp(-\Delta S_{\tau_+})} \\
&= \frac{E_{\theta_1}(1 - R(\infty)) \exp(-\Delta R(\infty))}{\Delta E_{\theta_1} \exp(-\Delta R(\infty))} \\
&= \frac{1}{\Delta} + \frac{1}{\Delta} \frac{\partial}{\partial b} \rho(\theta_1, \Delta).
\end{aligned}$$

Thus, since

$$\frac{\partial}{\partial b} \rho(\theta, b) = \sum_{m=0}^{\infty} (-1)^m \kappa_{m+1}(\theta) \frac{b^m}{m!},$$

it follows that Theorem 4 can be applied directly to provide the full asymptotic expansion for $E_{\theta_0} M$. Indeed, our analysis in Sections 3 to 5 yield an asymptotic expansion for $\kappa_n(\cdot)$ around zero, which in turn implies the expansion

$$E_{\theta_0} M = \frac{1}{\Delta} + \sum_{m=0}^{n} \sum_{j=0}^{n-m} (-1)^m \kappa_{m+1}^{(j)}(0) \frac{\theta_1(\Delta)^j}{j!} \frac{\Delta^m}{m!} + O(\Delta^{n+1})$$

valid for all $n \geq 0$. The explicit computation of the derivatives

$$\kappa_{m+1}^{(j)}(0) \triangleq \frac{d^j}{d\theta^j} \kappa_{m+1}(\theta) \bigg|_{\theta=0},$$

for $j, m \geq 0$, is discussed in Section 6.2.

Finally, the analytic extension of $\kappa_n(\cdot)$ and $r(\cdot)$ is a consequence of the following result.

PROPOSITION 1. *If $X_1$ has exponential moments and strongly nonlattice distribution, then, $I(\cdot)$ [defined as in (17) on a domain containing $[0, \upsilon) \times [0, \upsilon)$ with $\upsilon > 0$] can be analytically extended throughout a disc containing the origin in $\mathbb{C} \times \mathbb{C}$.*

Moreover, with the aid of Theorem 1, it follows easily [from (18) and the implicit function theorem] that $\Delta E_{\theta_0} M$ (initially defined for $\theta_0 < 0$) can be analytically extended [as a function of $\Delta(\theta_0)$] in a neighborhood of the origin in the complex plane.

**3. Short-time asymptotics for the Cauchy process.** The approach described in Section 2 suggests computing an asymptotic expansion for $r(\Delta)$ by developing appropriate expansions for $s(\Delta)$ and $I(\theta_1, \Delta)$. In this section we will show how asymptotics for $s(\Delta)$ can be obtained. Section 4 shows how asymptotics for $I(\theta, \Delta)$ [and, as a result, also for $I(\theta_1(\Delta), \Delta)$] can be reduced to the types of integrals considered here.

Since $s(b)$ is real for $b$ positive, it follows that the integral of the imaginary part of (16) must vanish. Hence, $s(b)$ equals the integral of the real part of (16), so that

$$(19) \quad \begin{aligned} s(b) &= \frac{1}{2\pi} \int_{-\infty}^{\infty} \frac{b}{b^2 + \lambda^2} \operatorname{Re} \log(2(1 - g(\lambda))\lambda^{-2}) \, d\lambda \\ &\quad - \frac{1}{2\pi} \int_{-\infty}^{\infty} \frac{b}{(b^2 + \lambda^2)\lambda} \operatorname{Im} \log(1 - g(\lambda)) \, d\lambda. \end{aligned}$$

Both of the above integrals take the form

$$(20) \quad \begin{aligned} K(b, f) &= \frac{1}{2\pi} \int_{-\infty}^{\infty} \frac{b}{b^2 + \lambda^2} f(\lambda) \, d\lambda \\ &= \frac{1}{2\pi} \int_{-\infty}^{\infty} \frac{1}{1 + \lambda^2} f(\lambda b) \, d\lambda. \end{aligned}$$

for suitably defined $f$. Note that if $Y = (Y(t) : t \geq 0)$ is a standard Cauchy process [so that $Y(1)$ is distributed as a standard Cauchy r.v.], $K(t, f)$ can then be represented as

$$K(t, f) = \tfrac{1}{2} E(f(Y(t)) | Y(0) = 0).$$

Hence, representing $K(t, f)$ as a power series in $t$ is equivalent to the development of short-time asymptotics of the Cauchy process. Such asymptotics are also of general analytical interest, because of their relevance to Fourier analysis. Integrals of the type (20) are closely related to "approximate identities of the Fejer type"; see page 31 of [5].

Let $\mathfrak{L}$ be the space of functions $f : \mathbb{R} \to \mathbb{C}$ for which $E|f(Y(1))|$ is finite and for which $f$ is infinitely differentiable at zero. For $f : \mathbb{R} \to \mathbb{C}$, let $\underline{f}$ be the symmetrization of $f$ defined via $\underline{f}(x) = (f(x) + f(-x))/2$. The following result provides our short-time asymptotic expansion for $K(t, f)$.

PROPOSITION 2. *Suppose $f$ belongs to $\mathfrak{L}$. Then, $K(\cdot, f)$ is infinitely differentiable at the origin and*

$$K^{(n)}(0, f) = \begin{cases} (-1)^{n/2} f^{(n)}(0), & n \text{ even,} \\ (-1)^{(n-1)/2} n! \dfrac{1}{2\pi} \displaystyle\int_{-\infty}^{\infty} (T_{(n-1)/2}\underline{f})(\lambda)\, d\lambda, & n \text{ odd,} \end{cases}$$

*where, for $j \geq 0$, $T_j$ acts on even functions in $\mathfrak{L}$ as*

$$T_j f(\lambda) = \frac{f(\lambda) - \sum_{k=0}^{2j-1} f^{(2k)}(0) \lambda^{2k}/(2k!)}{\lambda^{2j}}.$$

*Furthermore, the family of linear operators $(T_n : n \geq 0)$ is a commutative semigroup, so that $T_{n+m} = T_n T_m$ $m, n \geq 0$.*

REMARK. Note that the even derivatives of $f$ match those of $\underline{f}(\cdot)$. One might therefore be tempted to write the derivatives of $K(\cdot, f)$ in terms of integrals of $T_j f$ rather than $T_j \underline{f}$. The problem is that $T_j f$ typically has a singularity at the origin, unless the odd derivatives of $f$ at zero vanish. As a consequence, the integrals defining the derivative of $K(\cdot, f)$ may diverge if they were defined directly in terms of $f$. To avoid this, we use the symmetrization $\underline{f}$.

PROOF OF PROPOSITION 2. The fact that $T_n$ is a linear operator and forms a commutative semigroup is straightforward. To obtain the formula for the derivatives of $K(\cdot, f)$, at the origin, note that $K(\cdot, f) = K(\cdot, \underline{f})$ where $\underline{f}$ is the symmetrization of $f$ given by $\underline{f}(\cdot) = (f(\cdot) + f(-\cdot))/2$. Furthermore, if $f \in \mathfrak{L}$, then $\underline{f}$ is also in $\mathfrak{L}$. Observe that the dominated convergence theorem implies that

$$K(t, \underline{f}) \to \underline{f}(0)/2$$

as $t \searrow 0$. This motivates writing

$$K(t, \underline{f}) = \underline{f}(0)/2 + \frac{1}{2\pi} \int_{-\infty}^{\infty} \frac{t}{t^2 + \lambda^2} (\underline{f}(\lambda) - \underline{f}(0))\, d\lambda.$$

Since $E|\underline{f}(Y(1))|$ is finite, it follows that the above integrand is uniformly dominated by an integrable function for $|\lambda|$ bounded away from zero. On the other hand, $\underline{f}(\lambda) - \underline{f}(0) = O(\lambda^2)$ as $\lambda \to 0$, so the integrand is also uniformly (in $t$) dominated for $|\lambda|$ small. Hence, the dominated convergence theorem yields the conclusion that

$$K(t, f) = \underline{f}(0)/2 + \frac{t}{2\pi} \int_{-\infty}^{\infty} \frac{1}{t^2} (\underline{f}(\lambda) - \underline{f}(0))\, d\lambda + o(t)$$

as $t \to 0$. In fact,

$$
\begin{aligned}
K(t,f) &= \underline{f}(0)/2 + \frac{t}{2\pi} \int_{-\infty}^{\infty} \frac{1}{\lambda^2} (\underline{f}(\lambda) - \underline{f}(0))\, d\lambda \\
&\quad - \frac{t^2}{2\pi} \int_{-\infty}^{\infty} \frac{t}{t^2 + \lambda^2} \frac{(\underline{f}(\lambda) - \underline{f}(0))}{\lambda^2}\, d\lambda \\
&= \underline{f}(0)/2 + \frac{t}{2\pi} \int_{-\infty}^{\infty} (T_1 \underline{f})(\lambda)\, d\lambda - t^2 K(t, T_1 \underline{f}).
\end{aligned}
$$
(21)

If we apply (21) recursively to $K(\cdot, T_1\underline{f})$, $K(\cdot, T_2\underline{f}), \ldots$ we find that $K(t,f)$ satisfies

$$
K(t,f) = \sum_{j=0}^{n} (-1)^j \left( \frac{t^{2j}(T_j \underline{f})(0)}{2} + \frac{t^{2j+1}}{2\pi} \int_{-\infty}^{\infty} (T_{j+1}\underline{f})(\lambda)\, d\lambda \right)
$$
$$
+ (-1)^{n+1} t^{2(n+1)} K(t, T_{2(n+1)}\underline{f}),
$$

yielding the result. $\square$

With Proposition 2 in hand, our asymptotic expansion for $s(\Delta)$ follows immediately.

**4. Reducing the analysis to Cauchy process short-time asymptotics.** As we discussed earlier in Section 2, the backbone of our asymptotic analysis for $r(\Delta)$ is given by the relation $\rho(\theta, b) = s(b) + I(\theta, b)$. In Section 3 we studied how to develop asymptotics for $s(b)$. In this section we will study how to reduce the analysis of the remaining term $I(\theta, b)$ to that already studied in Section 3. Recall that

$$
I(\theta, b) = \frac{1}{2\pi} \int_{-\infty}^{\infty} \frac{-b}{(b + i\lambda)i\lambda} \log(1 - v(\theta, \lambda))\, d\lambda,
$$

where

$$
v(\theta, \lambda) = \frac{\lambda}{\lambda - 2i\phi'(\theta)} \left( 1 - \frac{2i\phi'(\theta)}{\lambda} - \frac{\gamma(\theta) - \gamma(\theta + i\lambda)}{1 - g(\lambda)} \right).
$$

A natural strategy that one could try to pursue now is to express the logarithm as a power series in $v(\theta, \lambda)$, followed by an expansion for $v$ as

(22)
$$
v(\theta, \lambda) = \sum_{n=0}^{\infty} v_n(i\lambda) \frac{\theta^n}{n!}.
$$

One could then apply Proposition 2 [as for (19)] to the real and imaginary parts in each of the resulting integrals that would appear as coefficients for $\theta^n$. However, the expansion (22) requires that the function $v$ be expressible as a joint power series in nonnegative powers of $\theta$ and $\lambda$. Unfortunately, the presence of the term $(\lambda - 2i\phi'(\theta))$ in the denominator of $v$ precludes the existence of such a joint power series.

To avoid this difficulty, we write $v$ as

$$v(\theta, \lambda) = \frac{\lambda H(\theta, \lambda)}{\lambda - 2i\phi'(\theta)},$$

so that

$$H(\theta, \lambda) = 1 - \frac{2i\phi'(\theta)}{\lambda} - \frac{\gamma(\theta) - \gamma(\theta + i\lambda)}{1 - g(\lambda)}.$$

The function $H(\cdot)$ is well behaved because the term $2i\phi'(\theta)/\lambda$ controls the behavior of $(\gamma(\theta) - \gamma(\theta + i\lambda))(1 - g(\lambda))^{-1}$ as $\lambda \searrow 0$. As a consequence, $H(\cdot)$ can be smoothly defined at $\lambda = 0$ via the relation $H(\theta, 0) = 1 - \phi''(\theta)$. Our next result describes the analytic structure of $H(\cdot)$.

PROPOSITION 3. *Let $D_{\eta/2} \triangleq \{z \in \mathbb{C} : |z| < \eta/2\}$ and, for $(z_1, z_2) \in D_{\eta/2}(D_{\eta/2} \cup \mathbb{R})$, put $H$*

$$H(z_1, z_2) = 1 - 2\frac{i\gamma'(z_1)}{z_2} - \frac{\gamma(z_1) - \gamma(z_1 + iz_2)}{1 - \gamma(iz_2)}.$$

*Then, for every $z_1 \in D_{\eta/2}$, the function $H(z_1, \cdot)$ is analytic on $D_{\eta/2}$. Similarly, for every $z_2 \in D_{\eta/2} \cup \mathbb{R}$, the function $H(\cdot, z_2)$ is analytic on $D_{\eta/2}$. Finally, $H(z_1, \lambda)$ can be represented as an absolutely and uniformly convergent series, for $\lambda \in \mathbb{R}$ and $z_1 \in D_{\eta/2}$, namely,*

$$(23) \qquad H(z_1, \lambda) = \sum_{k=1}^{\infty} h_k(i\lambda) \frac{z_1^k}{k!},$$

*where $h_k(i\lambda) \triangleq (\gamma^{(k)}(i\lambda) - \mu_k)/(1 - g(\lambda)) - (2i\mu_{k+1}/\lambda)$. In particular, this implies that*

$$\sup_{\lambda \in \mathbb{R}} |H(z_1, \lambda)| \to 0$$

*as $z_1 \to 0$.*

REMARK. Note that the function $\widetilde{H}(z_1, z_2) \triangleq H(z_1, z_2) - H(\theta, 0) = H(z_1, z_2) - 1 + \gamma''(z_1)$ satisfies the same properties stated for $H(\cdot)$ in Proposition 3 with $\widetilde{h}_k(i\lambda) \triangleq h_k(i\lambda) + \mu_{k+2}$, this follows from the analyticity of $\gamma(\cdot)$ and the fact that $\gamma''(0) = 1$. Moreover, observe that completely analogous analytic properties apply to the function $\widetilde{G}(z_1, z_2) = (\gamma''(z_1))^{-1}\widetilde{H}(z_1, z_2)$ defined on $D_{\eta/2} \times (D_{\eta/2} \cup \mathbb{R})$.

Note that $|\lambda/(\lambda - 2i\phi'(\theta))| = |\lambda|(\lambda^2 + (2\phi'(\theta))^2)^{-1/2} \leq 1$. It follows from Proposition 3 that, for $r > 0$ small enough,

$$\sup_{\theta \in (0,r)} \sup_{\lambda \in \mathbb{R}} |v(\theta, \lambda)| < 1.$$

Therefore, for all $0 < \theta < r$, we can proceed to expand $\log(1 - v)$ in powers of $v$ and formally integrate each term in the obtained expansion to express $I(\theta, b)$ in terms of integrals of the form

$$(24) \qquad J_k(a, b, f) = \frac{1}{2\pi} \int_{-\infty}^{\infty} \frac{-b}{(b + i\lambda)i\lambda} \left(\frac{i\lambda}{a + i\lambda}\right)^k f(i\lambda) \, d\lambda,$$

where $a, b > 0$, $f(i \cdot) \in \mathcal{L}$ and $k \geq 0$. Because $J_0(a, b, f) \triangleq J_0(b, f)$ can be written as

$$(25) \qquad \begin{aligned} J_0(b, f) &= \frac{1}{2\pi} \int_{-\infty}^{\infty} \frac{b}{b^2 + \lambda^2} (\operatorname{Re} f(i\lambda) - \lambda^{-1} \operatorname{Im} f(i\lambda)) \, d\lambda \\ &\quad + \frac{i}{2\pi} \int_{-\infty}^{\infty} \frac{b}{b^2 + \lambda^2} (\operatorname{Im} f(i\lambda) + \lambda^{-1} \operatorname{Re} f(i\lambda)) \, d\lambda, \end{aligned}$$

it follows that asymptotics for $J_0$ can be computed in terms of asymptotics for the $K$-type integrals that are subject of Proposition 2. In view of the development leading to (24), a key to our asymptotic expansion for $I(\theta, b)$ is therefore the reduction of integrals $J_k(a, b, f)$ for $k \geq 1$ to integrals of the form $J_0(b, f)$. A key identity in establishing this reduction step is the following.

LEMMA 1. *Suppose that $a, b \geq 0$. Then, for $m, n \geq 0$,*

$$\frac{1}{2\pi} \int_{-\infty}^{\infty} \frac{-b}{(b + i\lambda)i\lambda} \frac{(i\lambda)^{m+1}}{(a + i\lambda)^{m+n+1}} \, d\lambda = 0.$$

*Furthermore,*

$$\frac{1}{2\pi} \int_{-\infty}^{\infty} \frac{-1}{(1 + i\lambda)i\lambda} \log(1 + ai\lambda) \, d\lambda = 0.$$

PROOF. For $a, b > 0$, let the function of a complex variable $f(\cdot)$ be defined as

$$f(z) = \frac{-b}{(b + iz)iz} \frac{(iz)^{m+1}}{(a + iz)^{m+n+1}}.$$

Consider the contour (in the clockwise direction) $C(r) = C_1(r) + C_2(r)$, where $C_1(r) = \{re^{i\tau} : -\pi \leq \tau \leq 0\}$ and $C_2(r) = \{\lambda : \lambda \in [-r, r]\}$. Since $f$ is (complex) analytic on $\operatorname{Im}(z) \leq 0$, Cauchy's theorem yields

$$\frac{1}{2\pi} \int_{C(r)} \frac{-b}{(b + iz)iz} \frac{(iz)^{m+1}}{(a + iz)^{m+n+1}} \, dz = 0.$$

This, in turn, implies that

$$\frac{1}{2\pi}\int_{-r}^{r}\frac{-b}{(b+i\lambda)i\lambda}\frac{(i\lambda)^{m+1}}{(a+i\lambda)^{m+n+1}}\,d\lambda$$

$$=\frac{-1}{2\pi}\int_{C_1(r)}\frac{-b}{(b+iz)iz}\frac{(iz)^{m+1}}{(a+iz)^{m+n+1}}\,dz$$

$$=\frac{-1}{2\pi}\int_{-\pi}^{0}\frac{b(ir)^{m+1}e^{(m+1)\tau i}}{(b+ire^{\tau i})(a+ire^{\tau i})^{m+n+1}}\,d\tau.$$

Letting $r \to \infty$, we obtain (by virtue of dominated convergence) the first part of the lemma. For the second part, let us define

$$f_1(a) = \frac{1}{2\pi}\int_{-\infty}^{\infty} -((1+i\lambda)i\lambda)^{-1}\log(1+ai\lambda)\,d\lambda.$$

A routine dominated convergence argument, combined with our previous analysis, shows that

$$f_1'(a) = \frac{1}{2\pi}\int_{-\infty}^{\infty}\frac{-1}{(1+i\lambda)(1+ai\lambda)}\,d\lambda = 0.$$

The proof of the lemma is complete by observing that $f_1(a) \to 0$ as $a \searrow 0$. □

Let $\mathfrak{L}_0$ be the subspace of $\mathfrak{L}$ (recall the definition of $\mathfrak{L}$ preceding Proposition 2) for which $f(0) = 0$. Also, for $f \in \mathfrak{L}$, let $\widetilde{f}(\cdot) = f(\cdot) - f(0)$ ($\in \mathfrak{L}_0$). We are now ready to offer a proposition that reduces the evaluation of the integrals $J_k(a,b,f)$ for $k \geq 1$ to that of integrals such as $J_0(b,f)$, thereby permitting the application of the short-time asymptotics of Section 3.

PROPOSITION 4. *Suppose that $f \in \mathfrak{L}_0$. Then, for $k \geq 1$ and $n \geq 0$,*

$$(26) \qquad J_k(a,b,f) = J_0\left(b, \sum_{j=0}^{n}\binom{k+j-1}{j}(-a)^j \widetilde{T}_j \widetilde{f}\right) + bo(a^n),$$

*where the linear operator $\widetilde{T}_j$ ($j \geq 0$) acts on functions $\widetilde{f}(i\cdot) \in \mathfrak{L}_0$ as*

$$(\widetilde{T}_j \widetilde{f})(i\lambda) = \frac{\widetilde{f}(i\lambda) - \sum_{m=1}^{j}\widetilde{f}^{(m)}(0)(i\lambda)^m/m!}{(i\lambda)^j}.$$

*Moreover, the family of operators $(\widetilde{T}_j : j \geq 0)$ constitutes a commutative semigroup, so that $\widetilde{T}_m \widetilde{T}_n = \widetilde{T}_{m+n}$.*

REMARK. As for Proposition 2, one might be tempted to express the right-hand side of (26) in terms of $f$ rather that $\widetilde{f}$. However, $\widetilde{T}_j f$ is generally nonintegrable with respect to the kernel that defines $J_0$. Finally, note that, if all integrals are interpreted in terms of Cauchy principal value, one can apply Proposition 4 directly to functions that do not vanish at the origin by defining $J_0(b, f) = J_0(b, f(\cdot) - f(0)) + f(0)/2$.

PROOF OF PROPOSITION 4. That $(\widetilde{T}_j : j \geq 0)$ is a family of linear operators forming a commutative semigroup is immediate. By virtue of Lemma 1, it follows that

$$J_m(a, b, f) = \frac{1}{2\pi} \int_{-\infty}^{\infty} \frac{-b}{(b+i\lambda)i\lambda} \left(\frac{i\lambda}{a+i\lambda}\right)^m \widetilde{f}(i\lambda) \, d\lambda.$$

Observe that $\widetilde{f}(i\cdot)$ is now in the domain of the operators $\widetilde{T}_n$, $n \geq 1$. On the other hand, we can write

(27)
$$J_m(a, b, \widetilde{f}) = \frac{1}{2\pi} \int_{-\infty}^{\infty} \frac{-b}{(b+i\lambda)i\lambda} \widetilde{f}(i\lambda) \, d\lambda$$
$$+ \frac{1}{2\pi} \int_{-\infty}^{\infty} \frac{-b}{(b+i\lambda)i\lambda} \widetilde{f}(i\lambda) \left(\left(\frac{i\lambda}{a+i\lambda}\right)^m - 1\right) d\lambda.$$

Note that

$$\left(\frac{i\lambda}{a+i\lambda}\right)^m - 1 = -\sum_{k=1}^{m} \binom{m}{k} \frac{a^k (i\lambda)^{m-k}}{(a+i\lambda)^m}.$$

Once again, by appealing to Lemma 1 and to the definition of $\widetilde{T}_k \widetilde{f}$, it follows that, for $m \geq k \geq 1$,

$$a^k J_m(a, b, \widetilde{T}_k \widetilde{f}) \triangleq \frac{1}{2\pi} \int_{-\infty}^{\infty} \frac{-b}{(b+i\lambda)i\lambda} \frac{a^k (i\lambda)^{m-k}}{(a+i\lambda)^m} \widetilde{f}(i\lambda) \, d\lambda.$$

Combining this observation with (27), we obtain

(28)
$$J_m(a, b, f) = J_m(a, b, \widetilde{f})$$
$$= J_0(b, \widetilde{f}) - \sum_{k=1}^{m} \binom{m}{k} a^k J_m(a, b, \widetilde{T}_k \widetilde{f}).$$

The recursive relation (28) can now be expressed in operator form as

$$J_m(a, b, \widetilde{f}) = J_0(b, \widetilde{f}) + J_m(a, b, (1 - (1 + a\widetilde{T})^m)\widetilde{f}).$$

(Here, we have used the semigroup property of the family of operators $\widetilde{T}_m$.) Iterating the previous expression, we arrive at

$$J_m(a, b, f) = J_m(a, b, \widetilde{f})$$

$$
\begin{aligned}
(29) \quad &= \sum_{k=0}^{n} J_0(b,(1-(1+a\widetilde{T})^m)^k \widetilde{f}) \\
&\quad + J_m(a,b,(1-(1+a\widetilde{T})^m)^{n+1}\widetilde{f}) \\
&= J\left(b, \sum_{j=0}^{n} \binom{m+j-1}{j}(-a)^j \widetilde{T}_j \widetilde{f}\right) + bo(a^n),
\end{aligned}
$$

where the last equality in (29) has been obtained by using the semigroup property of the operators $\widetilde{T}_m$ and by noting that the coefficient of $a^j \widetilde{T}_j$ in (29) (for $j \leq n$) must match that of $x^j$ in the formal expansion of

$$
\begin{aligned}
p(x) &= \frac{1-(1-(1+x)^m)^{n+1}}{1-(1-(1+x)^m)} \\
&= \frac{1}{(1+x)^m} + O(x^{n+1}).
\end{aligned}
$$

That the error term in (29) is $bo(a^n)$ comes from the fact that $aJ_m(a,b,\widetilde{f}) = bo(1)$, as $a \searrow 0$, as it can be seen as follows:

$$
\begin{aligned}
|aJ_m(a,b,\widetilde{f})| &= \left|\frac{a}{2\pi}\int_{-\infty}^{\infty} \frac{-b(i\lambda)^{m-1}\widetilde{f}(i\lambda a)}{(b+i\lambda a)(1+i\lambda)^m}\,d\lambda\right| \\
&\leq \frac{b}{2\pi}\int_{-\infty}^{\infty}\left|\frac{\widetilde{f}(i\lambda a)}{\lambda(1+i\lambda)}\right|d\lambda = bo(1),
\end{aligned}
$$

where the last step follows by a dominated convergence argument. This concludes the proof of Proposition 2. □

Proposition 4, combined with our development for $K(t,\cdot)$ in Section 3, provides all the elements required to develop asymptotic expansions for integrals of the form $J_m(a,b,f)$. Since, as discussed earlier at a formal level, $I(\theta,b)$ can be expressed as a sum of terms such as $J_m(a,b,f)$, it follows that the whole asymptotic analysis of $r(\Delta)$ and $\rho(\theta,b)$ can be reduced to that of Section 3. A complete rigorous justification for this representation for $I(\theta,b)$ is one of the main issues discussed in Section 5.

**5. An asymptotic expansion for $I(\theta,b)$.** In Sections 3 and 4 we have developed the tools required to obtain asymptotic expansions, in powers of $b$, for $s(b)$ and $I(\theta,b)$. We have done this by showing that the problem can be reduced to short-time asymptotics for the Cauchy process. The purpose of this section is to make rigorous the expansion for $I(\theta,b)$, in powers of $\theta$, that was outlined in Section 4.

Noting the important role that functions vanishing at the origin plays in Proposition 4, it seems appropriate to define

(30)
$$\widetilde{H}(\theta,\lambda) \triangleq H(\theta,\lambda) - H(\theta,0) = H(\theta,\lambda) - 1 + \phi''(\theta)$$
$$= \sum_{k=1}^{\infty} \widetilde{h}_k(i\lambda)\frac{\theta^k}{k!},$$

where $\widetilde{h}_k(i\lambda) \triangleq (\gamma^{(k)}(i\lambda) - \mu_k)/(1 - g(\lambda)) - (2i\mu_{k+1}/\lambda) + \mu_{k+2}$ is such that $\widetilde{h}_k(0) = 0$. The next proposition shows how a simplified expression for $I(\theta,b)$ in terms of $\widetilde{H}$ can be obtained.

PROPOSITION 5. *Define* $\Psi(\theta) = 2\phi'(\theta)/\phi''(\theta)$. *Then,*

(31) $$I(\theta,b) = \frac{1}{2\pi}\int_{-\infty}^{\infty} \frac{-b}{(b+i\lambda)i\lambda} \log\left(1 - \frac{\phi''(\theta)^{-1}\lambda\widetilde{H}(\theta,\lambda)}{(\lambda - i\Psi(\theta))}\right) d\lambda.$$

PROOF. Just note that

$$\log(1 - v(\theta,\lambda)) = \log\left(1 - \frac{\lambda\widetilde{H}(\theta,\lambda)}{(\lambda - i\Psi(\theta)\phi''(\theta))} - \frac{\lambda(1 - \phi''(\theta))}{(\lambda - i\Psi(\theta)\phi''(\theta))}\right)$$
$$= \log\left(\frac{i\lambda/\Psi(\theta) + 1}{i\lambda\phi''(\theta)/\Psi(\theta) + 1}\right)$$
$$\quad + \log\left(1 - \frac{\phi''(\theta)^{-1}\lambda\widetilde{H}(\theta,\lambda)}{(\lambda - i\Psi(\theta))}\right).$$

Thus, (31) follows from Lemma 1 by noting that

$$\frac{1}{2\pi}\int_{-\infty}^{\infty} \frac{-b}{(b+i\lambda)i\lambda} \log\left(\frac{i\lambda/\Psi(\theta) + 1}{i\lambda\phi''(\theta)/\Psi(\theta) + 1}\right) d\lambda$$
$$= \frac{1}{2\pi}\int_{-\infty}^{\infty} \frac{-1}{(1+i\lambda)i\lambda} \log\left(\frac{i\lambda b/\Psi(\theta) + 1}{i\lambda b\phi''(\theta)/\Psi(\theta) + 1}\right) d\lambda = 0.$$

□

Additional simplifications reduce the complexity of the expansion for $I(\theta,b)$. In particular, the expression for the integral $J_0(b,f)$ simplifies when it is known that $J_0(b,f)$ is real; see (25). Fortunately, our analysis of $I(\theta,b)$ gives rise to such real-valued $J_0(b,f)$'s. To establish this result, we introduce the following family of functions.

DEFINITION. A function $f:\mathbb{R} \to \mathbb{C}$ is said to have the "parity property" if $\operatorname{Re} f(i\cdot)$ and $\operatorname{Im} f(i\cdot)$ are even and odd functions, respectively. The class of functions possessing the parity property will be denoted by $\mathcal{P}$.

Note that if $f(i\cdot)$ is in the domain of $J_0(b,\cdot)$ and $f$ possesses the parity property, then we must have that $\operatorname{Im} J_0(b,f)=0$ (since it corresponds to an integral on the real line of an odd integrable function). The family of functions enjoying the parity property has certain closure characteristics that will be useful for the rest of our development. These closure properties are discussed in the next proposition.

PROPOSITION 6. *The class $\mathcal{P}$ of functions forms an algebra on $\mathbb{R}$ (i.e., a vector space on $\mathbb{R}$ that is closed under product of functions). In addition, if $f \in \mathcal{P}$, then $1/f(\cdot)$ (defined on its domain of finiteness) also possesses the parity property. Finally, if $f$ is in the domain of $\widetilde{T}$ and has the parity property, then $\widetilde{T}f \in \mathcal{P}$.*

PROOF. Certainly $\mathcal{P}$ constitutes a vector space on $\mathbb{R}$ and it is almost immediate that $\widetilde{T}$ preserves the parity property. Now, if $f_1, f_2 \in \mathcal{P}$, then $\operatorname{Re}(f_1 f_2) = \operatorname{Re}(f_1)\operatorname{Re}(f_2) - \operatorname{Im}(f_1)\operatorname{Im}(f_2)$ must clearly be even. Similarly, $\operatorname{Im}(f_1 f_2)$ must be odd, which implies that $f_1 f_2 \in \mathcal{P}$. Finally, note that

$$\frac{1}{f} = \frac{\operatorname{Re}(f)}{\operatorname{Re}(f)^2 + \operatorname{Im}(f)^2} - i\frac{\operatorname{Im}(f)}{\operatorname{Re}(f)^2 + \operatorname{Im}(f)^2},$$

which immediately implies that $\operatorname{Re} 1/f$ and $\operatorname{Im} 1/f$ are even and odd functions respectively and, thus, $1/f \in \mathcal{P}$. □

We now present the main result of this section, which yields an expansion for $I(\theta, b)$ in powers of $\theta$ and coefficients involving only integrals of the form $J_0(b, f)$ with $f$ satisfying the parity property.

PROPOSITION 7. *For $k, m \geq 1$, let the coefficient multiplying $\theta^k$ in the power series representation of $\widetilde{G}(\theta,\lambda)^m \triangleq (\phi''(\theta)^{-1}\widetilde{H}(\theta,\lambda))^m$ be defined as $\widetilde{g}_{k,m}(i\lambda)$. Then, $\widetilde{g}_{k,m}(\cdot) \in \mathcal{P}$ can be recursively computed via*

$$\widetilde{g}_{k,m+1}(i\lambda) = \sum_{n=0}^{k} \widetilde{g}_{n+1,m}(i\lambda)\widetilde{g}_{k-n,1}(i\lambda).$$

*Consider $b > 0$ and let $\chi(\theta) = -\Psi(\theta)/\theta$. Then,*

$$(32) \qquad I(\theta, b) = \sum_{m=1}^{n} \theta^m \sum_{j=0}^{m-1} \chi(\theta)^j J_0(b, E_{j,m}) + bo(\theta^n),$$

*where $E_{j,m} \triangleq E_{j,m}(i\lambda)$, defined for $0 \leq j \leq m-1$ and $m \geq 1$ as*

$$E_{j,m} = -\sum_{k=0}^{m-j-1} \frac{1}{m-j-k} \binom{m-k-1}{j} \widetilde{T}_j \widetilde{g}_{k,m-j-k},$$

*satisfies the parity property.*

PROOF. Since $\gamma(i\lambda) = E_0 \cos(\lambda X_1) + iE_0 \sin(\lambda X_1)$, it follows that $\gamma^{(k)}(i\lambda) - \mu_k \in \mathcal{P}$, as does the function $1 - \gamma(i\lambda)$. By the closure properties described in Proposition 6, we may easily conclude that $\widetilde{g}_{k,1} \in \mathcal{P}$. A second application of Proposition 6 shows that $\widetilde{g}_{k,m} \in \mathcal{P}$ and $E_{j,m} \in \mathcal{P}$. The recursive expression provided for $\widetilde{g}_{k,m}$ follows from standard convolution operations of power series. For $n \geq 1$, define

$$\widetilde{G}_n(\theta, \lambda) \triangleq \sum_{k=1}^{n} \widetilde{g}_{k,1}(i\lambda) \frac{\theta^k}{k!}$$

and

$$I_n(\theta, b) \triangleq \frac{1}{2\pi} \int_{-\infty}^{\infty} \frac{-b}{(b+i\lambda)i\lambda} \log\left(1 - \frac{\widetilde{G}_n(\theta, \lambda)\lambda}{\lambda - i2\phi'(\theta)}\right) d\lambda.$$

Note that

$$\log\left(1 - \frac{\widetilde{G}(\theta, \lambda)\lambda}{\lambda - i2\phi'(\theta)}\right) - \log\left(1 - \frac{\widetilde{G}_n(\theta, \lambda)\lambda}{\lambda - i2\phi'(\theta)}\right)$$

$$= \log\left(1 - \frac{\lambda}{\lambda - i2\phi'(\theta)} \frac{(\widetilde{G}_n(\theta, \lambda) - \widetilde{G}(\theta, \lambda))}{(1 - \widetilde{G}_n(\theta, \lambda)\lambda(\lambda - i2\phi'(\theta))^{-1})}\right).$$

On the other hand, from the remark following Proposition 3 and because $\log(1+z) = z(1+\varepsilon(z))$ for $z \in \mathbb{C}$, where $|\varepsilon(z)| \leq |z|$ for $|z| \leq 1/2$ (see Proposition 8.46, [4]), we can see that there exists a constant $B > 0$ such that

$$|I(\theta, b) - I_n(\theta, b)|$$

$$\leq \frac{B}{2\pi} \int_{-\infty}^{\infty} \frac{b|\widetilde{G}_n(\theta, \lambda) - \widetilde{G}(\theta, \lambda)|}{(b^2 + \lambda^2)^{1/2}(\lambda^2 + (2\phi'(\theta))^2)^{1/2}} d\lambda.$$

Essentially, by making the change of variables $u = \lambda\theta$, we then see that, for all $\theta \in (0, \delta)$ for some $\delta > 0$, we have

$$|I(\theta, b) - I_n(\theta, b)| \leq \frac{Bb\theta^n}{2\pi} \int_{-\infty}^{\infty} \frac{|\widetilde{G}_n(\theta, \lambda\theta) - \widetilde{G}(\theta, \lambda\theta)|}{\theta^{n+1}|\lambda|(\lambda^2 + 1)^{1/2}} d\lambda.$$

It follows easily from the previous inequality and the dominated convergence theorem that

$$I(\theta, b) - I_n(\theta, b) = bo(\theta^n).$$

Using the expansion of $\log(1 + z)$ at $z = 0$ and a similar dominated convergence argument, we can write

$$I_n(\theta, b) = \frac{-1}{2\pi} \int_{-\infty}^{\infty} \frac{-b}{(b+i\lambda)i\lambda} \sum_{m=1}^{n} \frac{1}{m}\left(\frac{i\lambda}{i\lambda + \Psi(\theta)}\right)^m \widetilde{G}_n(\theta, \lambda)^m d\theta$$

$$+ bo(\theta^n)$$

(33)
$$= \frac{-1}{2\pi} \int_{-\infty}^{\infty} \frac{-b}{(b+i\lambda)i\lambda} \sum_{m=1}^{n} \frac{1}{m}\left(\frac{i\lambda}{i\lambda+\Psi(\theta)}\right)^m \sum_{k=0}^{n-m} \theta^{k+m}\widetilde{g}_{k,m}(i\lambda)\,d\lambda$$
$$+ bo(\theta^n).$$

Using Proposition 4 and (33), we obtain that

$I(\theta, b)$

(34)
$$= -\sum_{m=1}^{n} \theta^m \sum_{k=0}^{m-1} J_{m-k}\left(\Psi(\theta), b, \frac{\widetilde{g}_{k,m}}{m-k}\right)$$
$$+ bo(\theta^n)$$
$$= -\sum_{m=1}^{n} \theta^m \sum_{k=0}^{m-1} J_0\left(b, \sum_{j=0}^{n}\binom{m-k+j-1}{j}(-\Psi(\theta))^j \widetilde{T}_j \frac{\widetilde{g}_{k,m}}{m-k}\right)$$
$$+ bo(\theta^n)$$
$$= -\sum_{m=1}^{n} \theta^m \sum_{j=0}^{m-1} \theta^j \chi(\theta)^j J_0\left(b, \sum_{k=0}^{m-1}\binom{m-k+j-1}{j}\widetilde{T}_j \frac{\widetilde{g}_{k,m}}{m-k}\right)$$
$$+ bo(\theta^n)$$
$$= \sum_{m=1}^{n} \theta^m \sum_{j=0}^{m-1} \chi(\theta)^j J_0\left(b, -\sum_{k=0}^{m-j-1}\binom{m-k-1}{j}\widetilde{T}_j \frac{\widetilde{g}_{k,m}(i\lambda)}{m-j-k}\right)$$
$$+ bo(\theta^n),$$

which yields the desired conclusion. □

In view of the previous result, an explicit expression for the coefficients in the expansion for $J_0(\cdot, f)$, when $f$ satisfies the parity property, deserves special attention. Providing such explicit expressions is the aim of the next proposition.

PROPOSITION 8. *Suppose that $f(i\cdot) \in \mathcal{L}_0$ has the parity property. Then, $J_0(\cdot, f)$ is infinitely differentiable at zero and*

(35)
$$J_0^{(n)}(0, f) = \begin{cases} (-1)^{n/2}\left(f_{\mathrm{RE}}^{(n)}(0) \right.\\ \qquad\left. - \frac{n!}{2\pi}\int_{-\infty}^{\infty}(T_{n/2+1}f_{\mathrm{IM}})(\lambda)\,d\lambda\right), & n \text{ even,} \\ (-1)^{(n+1)/2}\left(\frac{f_{\mathrm{IM}}^{(n+1)}(0)}{(n+1)} \right.\\ \qquad\left. - \frac{n!}{2\pi}\int_{-\infty}^{\infty}(T_{(n+1)/2}f_{\mathrm{RE}})(\lambda)\,d\lambda\right), & n \text{ odd,} \end{cases}$$

where $f_{\mathrm{IM}}(i\lambda) = \operatorname{Im} f(i\lambda)\lambda^{-1}$ and $f_{\mathrm{RE}}(i\lambda) = \operatorname{Re} f(i\lambda)$.

PROOF. The proof follows by a direct application of Proposition 2 combined with the fact that $\operatorname{Re} J_0(b, f) = 0$. □

We close this section with some remarks that clarify how the expansion just derived for $I(\theta, b)$ can alternatively be viewed through the prism of a formal operator expansion. The analytic properties stated in Proposition 1 provide rigorous justification for the expansions outlined next. First, we note that if $\theta > 0$ is small enough and $b > 0$, we can formally write

$$(36) \qquad I(\theta, b) = -\sum_{k=1}^{\infty} \frac{1}{k} J_k(\Psi(\theta), b, \phi''(\theta)^{-k} \widetilde{H}^k(\theta, \cdot)).$$

Formally interpreting $(1 + a\widetilde{T})^{-m}$ as

$$(1 + a\widetilde{T})^{-m} = \sum_{k=0}^{\infty} \binom{m+k-1}{k} (-a)^k \widetilde{T}^k,$$

in combination with the expansion (26) developed for $J_k(a, b, f)$ and equality (36), allows us to write

$$I(\theta, b) = -\sum_{k=1}^{\infty} \frac{1}{k} J_0(b, \phi''(\theta)^{-k}(1 + \Psi(\theta)\widetilde{T})^{-k} \widetilde{H}^k(\theta, \cdot)).$$

If we introduce the convention that, for commutative operators $B_1(\theta)$, $B_2(\theta)$ and functions $F_1(\theta, \cdot)$, $F_2(\theta, \cdot)$, expressions of the form $B_1(\theta)F_1(\theta, \cdot) \times B_2(\theta)F_2(\theta, \cdot)$ (or any permutation of this form) are always interpreted as

$$(B_1(\theta)B_2(\theta))(F_1(\theta, \cdot)F_2(\theta, \cdot)),$$

then we can write

$$(37) \qquad I(\theta, b) = J_0(b, \log(1 - \phi''(\theta)^{-1}(1 + \Psi(\theta)\widetilde{T})^{-1} \widetilde{H}(\theta, \cdot))).$$

Expression (37) provides a convenient shorthand notation for the expansion of $I(\theta, b)$, in powers of $\theta$ and with coefficients in terms of integrals of the form $J_0(b, \cdot)$. In addition, note that, in order to recover the coefficients in the expansion for $I(\theta_1(\cdot), \cdot)$, one can apply formal differentiation to (37) in both arguments $\theta$ and $b$ [always having in mind that (37) is just a formalism representing a certain asymptotic expansion]. Hence, for example, one can obtain the first term in the expansion for $I(\theta_1(\cdot), \cdot)$ as

$$\partial_\Delta I(\theta_1(\Delta), \Delta)|_{\Delta=0} = \partial_\theta I(0, 0)\partial_\Delta \theta_1(0) + \partial_b I(0, 0),$$

22J. BLANCHET AND P. GLYNNwhere the formal derivatives applied to (37) must be interpreted using the formal operator convention introduced earlier. Thus, for example, if $B(\theta)$ is an operator of the form

$$B(\theta) = \sum_{k=0}^{\infty} b_k \theta^k \frac{\widetilde{T}^k}{k},$$

applied to a function $F(\theta, \lambda) = \sum f_k(i\lambda)\theta^k/k!$, we interpret the formal derivative $\partial_\theta \log(1 - B(\theta)F(\theta, \cdot))$ as

$$\partial_\theta \log(1 - B(\theta)F(\theta, \cdot)) = -\partial_\theta B(\theta)(1 - B(\theta)F(\theta, \cdot))^{-1} F(\theta, \cdot)$$
$$- B(\theta)(1 - B(\theta)F(\theta, \cdot))^{-1} \partial_\theta F(\theta, \cdot),$$

where

$$\partial_\theta B(\theta)(1 - B(\theta)F(\theta, \cdot))^{-1} F(\theta, \cdot)$$
$$= \sum_{k=0}^{\infty} (\partial_\theta B(\theta) B(\theta)^k) F(\theta, \cdot)^{k+1}$$

and, similarly,

$$B(\theta)(1 - B(\theta)F(\theta, \cdot))^{-1} \partial_\theta F(\theta, \cdot)$$
$$= \sum_{k=0}^{\infty} B(\theta)^{k+1} (F(\theta, \cdot)^k \partial_\theta F(\theta, \cdot)).$$

Thus, it is possible to combine this formalism with the expansion

$$J_0(b, f) = \sum_{n=1}^{m} J^{(n)}(0, f) b^n/n! + O(b^{m+1})$$

to recover the coefficients in the expansion for $I(\theta_1(\cdot), \cdot)$ in powers of $\Delta$.

**6. Expansions for $r(\Delta)$ and $E_\theta R^k(\infty)$.** In previous sections we developed all the elements required to rigorously compute a full asymptotic expansion for $r(\cdot)$ in powers of $\Delta$. In the first part of this section, as a summary, we indicate how the developments obtained in the previous three sections can be applied to provide an asymptotic expansion for $r(\cdot)$ in powers of $\Delta$. In view of the level of complexity in the computation of the constants $\beta_n$, the description in this section is intended to provide guidance for an easy-to-design practical implementation in a computational package such as Mathematica or Matlab. An efficient implementation of the procedure will appear elsewhere. In the second part of this section, also as a direct consequence of the analysis in the previous sections, we will develop a rigorous asymptotic expansion for the cumulants of $R(\infty)$ under $P_\theta$ in powers of $\theta$.

6.1. *The expansion for $r(\Delta)$.* An algorithm for computing $\beta_k$ for $k \leq n$ proceeds as follows:

1. Expand $s(\Delta)$ up to terms of order $O(\Delta^{n+1})$ using Proposition 8.
2. Similarly, expand the functions $J_0(\cdot, E_{j,m})$ up to terms $O(\Delta^{n-m})$ with $0 \leq j \leq m-1$ and $1 \leq m \leq n$. This also can be done by applying Proposition 8, since $E_{j,m}$ has the parity property.
3. Finally, the terms obtained can be combined with an expansion for $\theta_1(\Delta)$ up to terms of order $O(\Delta^{n+1})$. Such an expansion can be easily obtained using the implicit function theorem and therefore is omitted.

Observe that the previous algorithm provides an asymptotic expansion for $r(\cdot)$ in powers of $\Delta$. However, because of Theorem 1, we actually have that this asymptotic expansion converges absolutely in a neighborhood of the origin.

As a simple application of the previous expansion, we show that $\beta_2 = 0$.

PROPOSITION 9. *Suppose that $X_1$ has exponential moments and is strongly nonlattice. Then*
$$r(\Delta) = -\Delta \beta_1 + O(\Delta^3).$$

PROOF. We only need to show that $\beta_2 = 0$. Note that, by virtue of Proposition 8, the coefficient multiplying $\Delta^2$ in the expansion of $s(\Delta)$ equals
$$s_2 = \frac{1}{2\pi} \int_{-\infty}^{\infty} \frac{1}{\lambda^2} \left( \frac{\operatorname{Im} \log(1 - g(\lambda))}{\lambda} - \mu_3 \right)$$
$$- \left( \frac{\mu_4}{12} - \frac{\mu_3^2}{18} \right).$$

In order to show that $\beta_2 = 0$, it suffices to show that $\theta_1 J(\Delta, E_{0,1}) \sim -\Delta^2 s_2$ or [since $J(\Delta, E_{j,m}) = O(\Delta)$, $\theta_1/2 \sim \Delta$ and $\phi''(\theta_1) \sim 1$] that $\Delta J(\Delta, E_{0,1}) \sim -2\Delta^2 s_2$, where

$$\Delta J(\Delta, E_{0,1}) = \frac{1}{2\pi} \int_{-\infty}^{\infty} \frac{-\Delta}{(\Delta + i\lambda)i\lambda} \left( \frac{\gamma'(i\lambda)}{1 - g(\lambda)} - \frac{2i}{\lambda} + \mu_3 \right) d\lambda$$

$$(38) \qquad = \frac{1}{\pi} \int_0^{\infty} \frac{\Delta^2}{(\Delta^2 + \lambda^2)} \operatorname{Re}\left( \frac{\gamma'(i\lambda)}{1 - g(\lambda)} - \frac{2i}{\lambda} + \mu_3 \right) d\lambda$$

$$(39) \qquad - \frac{\Delta}{\pi} \int_0^{\infty} \frac{\Delta^2}{(\Delta^2 + \lambda^2)\lambda} \operatorname{Im}\left( \frac{\gamma'(i\lambda)}{1 - g(\lambda)} - \frac{2i}{\lambda} + \mu_3 \right) d\lambda.$$

Note that $g'(\lambda) = i\gamma'(i\lambda)$ and that $\operatorname{Im} \log(2\lambda^{-2}) = 0$; hence, we can write
$$\operatorname{Re}\left( \frac{\gamma'(i\lambda)}{1 - g(\lambda)} - \frac{2i}{\lambda} + \mu_3 \right) = -\operatorname{Im} \frac{d}{d\lambda}(\log(2(1 - g(\lambda))\lambda^{-2}) - \mu_3 i\lambda),$$

which implies, using integration by parts, that the integral in (38) equals

$$
\begin{aligned}
&\frac{-1}{\pi}\int_0^\infty \frac{2\lambda\Delta^2}{(\Delta^2+\lambda^2)^2}\operatorname{Im}(\log(2(1-g(\lambda))\lambda^{-2})-\mu_3 i\lambda)\,d\lambda\\
&\sim -\frac{\Delta^2}{\pi}\int_0^\infty \frac{2}{\lambda^2}\left(\frac{\operatorname{Im}\log(1-g(\lambda))}{\lambda}-\mu_3\right)d\lambda,
\end{aligned}
\tag{40}
$$

where (40) has been obtained using dominated convergence and simple manipulations. It follows from Proposition 2 and a first-order asymptotic expansion of $E_{0,1}(i\lambda)$ that (39) equals $-\Delta^2(\mu_3^2/9-\mu_4/6)$. Combining this last estimate together with (40) into (38) and (39) yields $\Delta J(\Delta, E_{0,1})\sim -2\Delta^2 s_2$, which is exactly what we wanted to show to conclude that $\beta_2=0$. □

6.2. *The expansion for* $E_\theta R(\infty)^k$ *as* $\theta\searrow 0$. We shall provide asymptotics for $E_\theta R(\infty)^k = E_\theta(S_{\tau_+}^k)/(k!E_\theta(S_{\tau_+}))$ via the cumulants $(\kappa_j(\theta):j\geq k)$ of $R(\infty)$ under $P_\theta$. In particular, these estimates yield the proof of Theorem 4 stated in Section 2. The idea is to develop an asymptotic expansion, in powers of $b$, for $s(b)$ and $I(\theta,b)$ respectively and to match coefficients in the expression

$$
\begin{aligned}
\rho(\theta,b) &= -\kappa_1(\theta)b + \kappa_2(\theta)b^2/2 - \kappa_3(\theta)b^3/3! + \cdots\\
&= s(b) + I(\theta,b).
\end{aligned}
\tag{41}
$$

In order to perform this task, we will take advantage of Proposition 7 as follows; first let us define, for $k\geq 1$, $\alpha_{k,j,m}=J_0^{(k)}(0,E_{j,m})/k!$ (which can be explicitly computed via Proposition 8). With this notation, we can write, for $l,n\geq 1$,

$$
\begin{aligned}
I(\theta,b) &= \sum_{m=1}^n \theta^m \sum_{j=0}^{m-1}\chi(\theta)^j\left(\sum_{k=1}^l \alpha_{k,j,m}b^k + O(b^{l+1})\right) + bo(\theta^n)\\
&= \sum_{k=1}^l b^k \sum_{m=1}^n \sum_{j=0}^{m-1} \theta^m \chi(\theta)^j \alpha_{k,j,m} + \theta O(b^{l+1}) + bo(\theta^n).
\end{aligned}
$$

Therefore, we obtain that, for all $s,n\geq 1$, $\kappa_s(\theta)$ satisfies

$$
\kappa_s(\theta) = (-1)^s\left(\kappa_s(0) + s!\sum_{m=1}^n \sum_{j=1}^{m-1}\theta^m\chi(\theta)^j\alpha_{s,j,m}\right) + O(\theta^{n+1}).
$$

Consequently, $\kappa_n(\cdot)$ is an infinitely differentiable function at $\theta=0$ and for $m\geq 0$ and $n\geq 1$, we have

$$
\frac{\kappa_n^{(m)}(0)}{n!} = (-1)^n\frac{\kappa_n(0)}{n!} + \sum_{s=0}^{m-1}\sum_{j=0}^{m-1-s}\chi_{s,j}\alpha_{n,j,m-s},
$$

where, for $n, j \geq 1$, $\chi_{n,j}$ is the coefficient multiplying $\theta^n$ in the expansion for $\chi(\theta)^j$. In particular, the $\chi_{n,j}$ can be computed recursively as

$$\chi_{n,j+1} = \sum_{n=0}^{k} \chi_{n,j} \chi^{(k-n)}(0)/(k-n)!,$$

with $\chi_{n,1} = \chi^{(n)}(0)/n!$.

## 7. Technical proofs.

PROOF OF THEOREM 3. Using Lemma 1, we can add

$$0 = \frac{1}{2\pi} \int_{-\infty}^{\infty} \frac{b}{(b+i\lambda)i\lambda} \log(1 + i\lambda/2\phi'(\theta)) \, d\lambda$$

to expression (13) for $\rho(\theta, b)$ to obtain

$$\rho(\theta, b) = \frac{1}{2\pi} \int_{-\infty}^{\infty} \frac{-b}{(b+i\lambda)i\lambda} \log\left(\frac{\gamma(\theta) - \gamma(\theta + i\lambda)}{-i\phi'(\theta)\lambda(1 + i\lambda/2\phi'(\theta))}\right) d\lambda$$

$$= \frac{1}{2\pi} \int_{-\infty}^{\infty} \frac{-b}{(b+i\lambda)i\lambda} \log\left(\frac{2(\gamma(\theta) - \gamma(\theta + i\lambda))}{\lambda(\lambda - 2i\phi'(\theta))}\right) d\lambda,$$

yielding the conclusion of the theorem. □

PROOF OF PROPOSITION 3. It follows immediately, by a Taylor series expansion of $\gamma(\cdot)$, that a series representation for $H$ can be written (for fixed $\lambda$ and $\theta$ such that $0 < |\lambda| + |\theta| < \eta$) as

$$H(\theta, \lambda) = 1 - \frac{2i\phi'(\theta)}{\lambda} - \frac{\gamma(\theta) - \gamma(\theta + i\lambda)}{1 - g(\lambda)}$$

$$= 1 - \frac{2i}{\lambda} \sum_{k=1}^{\infty} \mu_{k+1} \frac{\theta^k}{k!} - \frac{1}{1-g(\lambda)} \sum_{k=0}^{\infty} (\mu_k - \gamma^{(k)}(i\lambda)) \frac{\theta^k}{k!}$$

$$= \sum_{k=1}^{\infty} h_k(i\lambda) \frac{\theta^k}{k!}.$$

In fact, the functions $h_k(i\cdot)$ can be analytically extended throughout the disc $D_{\eta/2} = \{z \in \mathbb{C} : |z| < \eta/2\}$. This is easily seen as follows, recall that $\gamma(\cdot)$ [and therefore $\gamma^{(k)}(\cdot)$] are analytic on $\mathcal{N}$ (defined in Section 2). Also, observe that $1 - \gamma(iz) \sim z^2/2$ and $\gamma^{(k)}(iz) - \mu_k \sim iz\mu_{k+1}$ as $z \to 0$. Thus, $(\gamma^{(k)}(iz) - \mu_k)/(1 - \gamma(iz))$ possesses a simple pole at 0 with residue equal to $2i\mu_{k+1}$, which implies that the natural extension of $h_k$ defined as

$$h_k(iz) = \frac{\gamma^{(k)}(iz) - \mu_k}{1 - \gamma(iz)} - \frac{2i\mu_{k+1}}{z}$$

$$= \frac{(\gamma^{(k)}(iz) - \mu_k)z - 2i\mu_{k+1}(1 - \gamma(iz))}{(1 - \gamma(iz))z}$$

is analytic on $D_{\eta/2}$. Now, by virtue of the maximum principle (see, e.g., [14], page 253), we have that if $\delta > 0$ is suitably small,

$$\sup_{|z| \leq \delta} |h_k(iz)| \leq \sup_{|z|=\delta} |h_k(iz)|.$$

Since $\gamma(z)$ is a nonconstant analytic function defined on $D_{\eta/2}$ (which is an open set and thus has an accumulation point), then $1 - \gamma(z)$ has an isolated zero at $z = 0$. Thus, it is possible to choose $\delta > 0$ in such a way that

$$\inf_{|z|=\delta} |1 - \gamma(iz)| > \varepsilon > 0,$$

for some $\varepsilon > 0$. Consequently,

$$\sup_{|z| \leq \delta} |h_k(iz)| \leq \sup_{|z|=\delta} |h_k(iz)|$$

$$\leq \frac{1}{\varepsilon \delta} \sup_{|z|=\delta} |(\gamma^{(k)}(iz) - \mu_k)z + 2\mu_{k+1}(1 - \gamma(iz))|.$$

Observe that, for $|z| < \eta/2$, $\gamma^{(k)}(z) = E_0(X^k \exp(zX))$. Therefore, if $z = x + iy$, with $|z| = \delta$,

$$|\gamma^{(k)}(iz)| \leq E_0(|X|^k |\exp(izX)|)$$
$$= E_0(|X|^k |\exp(yX)|)$$
$$\leq E_0(|X|^k \exp(\delta |X|)).$$

A similar bound can be obtained for $\gamma(z)$ and we can conclude that $\exists B > 0$ such that

$$\sup_{|z| \leq \delta'} |h_k(iz)| \leq B(E_0(|X|^k(\exp(\delta |X|) + 1)) + E_0(|X|^{k+1})(1 + E_0 \exp(\delta |X|))).$$

Now, suppose that $\delta < \eta/2$. Then, if $z_1 \in D_{\eta/2}$, we can define

$$BE_0(|X|^{k+1})(1 + E_0 \exp(\delta |X|)) \frac{z_1^k}{k!} \triangleq N_1(z_1)$$

in such a way that the previous series converges absolutely and uniformly on $D_{\eta/2}$. Similarly, we can define

$$B \sum_{k=1}^{\infty} E_0(|X|^k(\exp(\delta |X|) + 1)) \frac{z_1^k}{k!} = BE_0\left(\sum_{k=1}^{\infty} |X|^k \frac{z_1^k}{k!}(\exp(\delta |X|) + 1)\right)$$
$$= BE_0((\exp(z_1 |X|) - 1)(\exp(\delta |X|) + 1))$$
$$\triangleq N_2(z_1).$$

Note that, for $j = 1, 2$, $N_j(z_1) \to 0$ as $z_1 \to 0$. On the other hand, since $g(\lambda)$ is strongly nonlattice, we have that

$$\sup_{|\lambda| \geq \delta} |h_k(i\lambda)| = \sup_{|\lambda| \geq \delta} \left| \frac{\gamma^{(k)}(i\lambda) - \mu_k}{1 - g(\lambda)} - \frac{2i\mu_{k+1}}{\lambda} \right|$$
$$\leq B(E_0(|X|^k) + E_0(|X|^{k+1})),$$

if $B < \infty$ is big enough. The previous estimates imply that there exist constants $0 < M_k \leq B(E_0(|X|^k(\exp(\delta|X|) + 1)) + E_0(|X|^{k+1})(1 + E_0 \exp(\delta|X|)))$ such that

$$|h_k(iz_2)| \leq M_k$$

for $z_2 \in \mathbb{R} \cup D_{\eta/2}$ and $|\sum_{k=1}^\infty M_k \frac{z_1^k}{k!}| \leq \sum_{k=1}^\infty |M_k \frac{z_1^k}{k!}| < \infty$ for $z_1 \in D_{\eta/2}$. Thus, using the Weierstrass M test, we obtain the validity of (23). Finally, the invoked Weierstrass M test, combined with the analytic functions convergence theorem (see Theorem 10.28, page 214, of [14]), yields the analyticity of $H(z_1, \cdot)$ on $\mathbb{R} \cup D_{\eta/2}$ (for $z_1 \in D_{\eta/2}$) and similarly for $H(\cdot, z_2)$ on $D_{\eta/2}$ (for $z_2 \in \mathbb{R} \cup D_{\eta/2}$). □

PROOF OF PROPOSITION 1. We start by writing

$$I(\theta, b) = \frac{1}{2\pi} \int_{-\infty}^{\infty} \frac{-b}{(b + i\lambda)i\lambda} \log\left(1 - \frac{H(\theta, \lambda)\lambda}{\lambda - 2\phi'(\theta)i}\right) d\lambda.$$

The strategy will be to study this integral on $\{|\lambda| < \delta\}$ and $\{|\lambda| \geq \delta\}$ separately (where $\delta > 0$ is some convenient small number to be characterized later):

$$(42) \quad I(\theta, b) = -\frac{1}{2\pi} \int_{-\delta}^{\delta} \frac{b}{(b + i\lambda)} \frac{1}{i\lambda} \log\left(1 - \frac{H(\theta, \lambda)\lambda}{\lambda - 2\phi'(\theta)i}\right) d\lambda$$

$$(43) \quad\quad\quad - \frac{1}{2\pi} \int_{|\lambda| \geq \delta} \frac{b}{(b + i\lambda)} \frac{1}{i\lambda} \log\left(1 - \frac{H(\theta, \lambda)\lambda}{\lambda - 2\phi'(\theta)i}\right) d\lambda.$$

Let us define $I_A(\theta, b)$ and $I_B(\theta, b)$ as (42) and (43) respectively. Suppose that $0 < b < \delta < \eta/2$. By making $u = b\lambda$, we can write

$$I_A(\theta, b) = -\frac{1}{2\pi} \int_{-\delta}^{\delta} \frac{b}{(b + i\lambda)} \frac{1}{i\lambda} \log\left(1 - \frac{H(\theta, \lambda)\lambda}{\lambda - 2\phi'(\theta)i}\right) d\lambda.$$

Let $C = \{w \in \mathbb{C} : |w| \leq \delta\} \cap \{\text{Im}(w) \leq 0\}$, and observe that, by virtue of Proposition 3, we can pick $\delta_1 > 0$ in such a way that, for all $0 < \theta < \delta_1$, the function

$$f_1(w) = \frac{b}{(b + iw)} \frac{1}{iw} \log\left(1 - \frac{H(\theta, w)w}{w - 2\phi'(\theta)i}\right)$$

is analytic on $C$. Thus, applying Cauchy's theorem to the contour enclosing $C$, we obtain

$$
\begin{aligned}
I_A(\theta, b) &= -\frac{1}{2\pi} \int_{-\pi}^{0} \frac{b}{(b + i\delta e^{i\lambda})} \frac{i\delta e^{i\lambda}}{i\delta e^{i\lambda}} \\
&\qquad \times \log\left(1 - \frac{H(\theta, \delta e^{i\lambda})\delta e^{i\lambda}}{\delta e^{i\lambda} - i}\right) d\lambda \\
&= \frac{1}{2\pi} \int_{-\pi}^{0} \frac{ib\delta^{-1}e^{-i\lambda}}{(1 - ib\delta^{-1}e^{-i\lambda})} \\
&\qquad \times \log\left(1 - \frac{H(\theta, \delta e^{i\lambda})}{1 - i2\phi'(\theta)\delta^{-1}e^{-i\lambda}}\right) d\lambda.
\end{aligned}
\tag{44}
$$

Equality (44) has been obtained by simple algebraic manipulations. Observe that the previous expression in combination with Proposition 3 and the analyticity of the functions $\phi'(\theta)$ ($\sim 0$) at zero immediately gives that $I_A(\theta, b)$ can be represented as an absolutely convergent double power series in $\theta$ and $b$ on the set $0 < |\theta| + |b| < \delta_2$ for some $\delta_2 > 0$. Indeed, if we pick $\delta_2$ small enough, it is possible to provide an explicit power series representation for $I_A(\theta, b)$ by using the expansion of $\log(1 - w)$ at $w = 0$ in combination with the series representation (23) for the function $H(\theta, \lambda)$ derived in Proposition 3 and a Taylor expansion of $(1 - w)^{-1}$ around $w = 0$.

The analysis of $I_B(\theta, b)$ is easier,

$$
I_B(\theta, b) = \frac{1}{2\pi} \int_{|\lambda| \geq \delta} \frac{b}{(1 - b\lambda^{-1})} \frac{1}{\lambda^2} \log\left(1 - \frac{H(\theta, \lambda)}{1 - i2\phi'(\theta)\lambda^{-1}}\right) d\lambda.
$$

Hence, in order to show that $I_B(\cdot)$ can be written as an absolutely convergence double power series in a neighborhood of the origin, it suffices to show (by Fubini's theorem) that

$$
\int_{|\lambda| \geq \delta} \sum_{k,j,m \geq 0}^{\infty} \binom{m+k}{k} \frac{b^{j+1}(2(E(\exp(\theta|X|) - 1 - |X|)))^m}{(k+1)|\lambda|^{j+2+m}} \\
\times \left(\sum_{s=1}^{\infty} |h_s(i\lambda)| \frac{\theta^s}{s!}\right)^{k+1} d\lambda
$$

is finite for all nonnegative $\theta$ and $b$ such that $\theta + b < \delta_3$ for some $\delta_3 > 0$. But this fact follows easily from Proposition 3, first note, by the change of variables $\lambda = u\delta$, that the previous expression equals

$$
\int_{|u| \geq 1} \sum_{k,j,m \geq 0}^{\infty} \binom{m+k}{k} \frac{b^{j+1}(2(E(\exp(\theta|X|) - 1 - |X|)))^m}{\delta^{j+m+1}(k+1)|u|^{j+2+m}}
$$

$$\times \left( \sum_{s=1}^{\infty} |h_s(i\lambda\delta)| \frac{\theta^s}{s!} \right)^{k+1} du,$$

now pick $\delta_3$ small enough so that $0 < \max(b, 2(E(\exp(\theta|X|) - 1 - |X|))) < \delta_3 < \delta$ (if $\theta + b < \delta_3$), and use Proposition 3 to conclude that one $\delta_3$ can be chosen so that $\sum_{s=1}^{\infty} |h_s(i\lambda\delta)| \frac{\theta^s}{s!} < c < 1 - \delta_3/\delta$. Therefore, we can bound the previous sum by

$$\int_{|u|\geq 1} \sum_{k,j,m\geq 0} \binom{m+k}{k} \frac{(\delta_3/\delta)^{j+m+1}}{(k+1)|u|^2} c^{k+1} \, du$$

$$\leq \frac{2}{3} \frac{1}{1-\delta_3/\delta} \left| \log\left(1 - \frac{c}{1-\delta_3/\delta}\right) \right| < \infty.$$

The conclusions obtained for both $I_A(\cdot)$ and $I_B(\cdot)$ indicate that, for all $0 \leq \theta, b \leq \upsilon$ (for some $\upsilon > 0$), $I(\theta, b)$ can be written as

$$I(\theta, b) = \sum_{j,k \geq 1} \theta^j b^k I_{jk},$$

where the previous series converges absolutely on the specified region on $\theta$ and $b$. The previous expression provides the natural analytic extension of $I(\cdot)$ on $D_\upsilon^2 = \{(z_1, z_2) \in \mathbb{C} \times \mathbb{C} : |z_1| + |z_2| < \upsilon\}$. □

PROOF OF THEOREM 1. Since

$$\exp(s(\Delta)) = \frac{1 - E_0(\exp(-\Delta S_{\tau_+}))}{\Delta E_0(S_{\tau_+})} = E_0(\exp(-\Delta R(\infty))),$$

the analytic extension of the term $s(\Delta)$ follows from that of the right-hand side, which comes from the fact that $S_{\tau_+}$ has exponential moments (see [2]). Thus, since $r(\Delta) = s(\Delta) + I(\theta_1(\Delta), \Delta)$, we just have to analyze $I(\theta_1(\Delta), \Delta)$. However, from the implicit function theorem, we know that $\theta_1(\cdot)$ is analytic in the neighborhood of the origin, thus, the analytic functions convergence theorem (see Theorem 10.28, page 214, of [14]) combined with Proposition 1 yields the desired conclusion. □

PROOF OF THEOREM 4. From Theorem 1, we know that, for $0 \leq \theta, b \leq \upsilon$ (for some $\upsilon > 0$),

$$I(\theta, b) = \sum_{j=1}^{\infty} b^j I_{\cdot,j}(\theta),$$

where each function $I_{\cdot,j}(\theta)$ can be expanded in absolutely convergent power series for $0 \leq \theta \leq \upsilon$, and thus can be analytically extended throughout a neighborhood of the origin in the complex plane. But,
$$\rho(\theta, b) = -\kappa_1(\theta)b + \kappa_2(\theta)b^2/2 - \kappa_3(\theta)b^3/3! + \cdots$$
$$= s(b) + I(\theta, b),$$
where $s(\cdot)$ is (real) analytic at zero. Hence, the conclusion of the theorem follows immediately by matching coefficients. $\square$

Next, we show that if the distribution of $X_1$ is symmetric, then, for $n \geq 1$, $\beta_{2n} = 0$.

PROOF OF THEOREM 2. As we discussed before, all that we need to show is that $\beta_{2n} = 0$. We have shown that an absolutely convergent power series representation is possible for $r(\Delta)$ when $\Delta$ is small, thus, it suffices to show that if $0 < \Delta < \delta$ (where $\delta > 0$ is suitably small), then an asymptotic expansion for $r(\Delta)$ is given in odd powers of $\Delta$ only. Using the integral expression (14), integrating on $|\lambda| \leq \delta$ and $|\lambda| > \delta$, we can write

$$(45) \quad r(\Delta) = \frac{1}{2\pi} \int_{|\lambda| < \delta} \frac{-\Delta}{(\Delta + i\lambda)i\lambda} \log\left(\frac{2(\gamma(\theta_1) - \gamma(\theta_1 + i\lambda))}{\lambda(\lambda - 2i\phi'(\theta_1))}\right) d\lambda$$

$$(46) \qquad + \frac{1}{2\pi} \int_{|\lambda| \geq \delta} \frac{-\Delta}{(\Delta + i\lambda)i\lambda} \log\left(\frac{2(\gamma(\theta_1) - \gamma(\theta_1 + i\lambda))}{\lambda(\lambda - 2i\phi'(\theta_1))}\right) d\lambda.$$

Define by $A(\Delta)$ and $B(\Delta)$ the integrals appearing in expressions (45) and (46), respectively. We first analyze $A(\Delta)$. Using a similar argument as in the proof of Theorem 1, we see that

$$A(\Delta) = \frac{1}{2\pi} \int_{C_1} \frac{\Delta}{(\Delta + iz)iz} \log\left(\frac{2(\gamma(\theta_1) - \gamma(\theta_1 + iz))}{z(z - 2i\phi'(\theta_1))}\right) dz,$$

where the trajectory $C_1$ is defined as $C_1 = \{\delta e^{i\lambda} : \lambda \in [0, -\pi)\}$. Also, define the trajectory $C_2 = \{\delta e^{i\lambda} : \lambda \in [-\pi, 0)\}$. The proof of the theorem will be complete if we show that $A(\Delta)$ is an odd function. That is, we must show that $A(\Delta) = -A(-\Delta)$. Note that

$$(47) \quad \begin{aligned} -A(-\Delta) &= \frac{-1}{2\pi} \int_{C_1} \frac{-\Delta}{(-\Delta + iz)iz} \log\left(\frac{2(\gamma(-\theta_1) - \gamma(-\theta_1 + iz))}{z(z - 2i\phi'(-\theta_1))}\right) dz \\ &= \frac{1}{2\pi} \int_{C_2} \frac{-\Delta}{(\Delta + iw)iw} \log\left(\frac{2(\gamma(\theta_1) + \gamma(\theta_1 + iw))}{w(w - 2i\phi'(\theta_1))}\right) dw. \end{aligned}$$

Equality (47) was obtained by making the change of variables $-w = z$ and using that $\gamma(\theta_1)$ and $\phi'(\theta_1)$ are even and odd functions of $\theta_1$, respectively. In view of (47), in order to show that $A(\Delta) = -A(-\Delta)$, it suffices to show that

$$0 = \frac{1}{2\pi} \int_C \frac{\Delta}{(\Delta + iw)iw} \log\left(\frac{2(\gamma(\theta_1) - \gamma(\theta_1 + iw))}{w(w - 2i\phi'(\theta_1))}\right) dw,$$

where $C = C_1 + C_2$ is the contour corresponding to the circle with radius $\delta$. Now,

$$\frac{1}{2\pi} \int_C \frac{\Delta}{(\Delta + iw)iw} \log\left(\frac{2(\gamma(\theta_1) - \gamma(\theta_1 + iw))}{w(w - 2i\phi'(\theta_1))}\right) dw$$

(48) $$= \frac{1}{2\pi} \int_{-C} \frac{\Delta}{w(w - i\Delta)} \log\left(\frac{2(\gamma(\theta_1) - \gamma(\theta_1 + iw))}{w(w - i\Delta)}\right) dw$$

(49) $$+ \frac{1}{2\pi} \int_{-C} \frac{\Delta}{w(w - i\Delta)} \log\left(\frac{w - i\Delta}{w - i2\phi'(\theta_1)}\right) dw.$$

We will show that both terms (48) and (49) vanish. We first consider (49). For $\gamma \in [0,1]$ and $a \in [-\delta, \delta]$, define $f(\gamma)$ as

$$f(\gamma) = \frac{1}{2\pi} \int_{-C} \frac{\Delta}{w(w - i\Delta)} \log(\gamma w - ia) \, dw.$$

Using residue calculus (see [14], page 224), it is easy to see that $f(0) = 0$. Indeed, a standard dominated convergence argument yields

$$f'(\gamma) = \frac{1}{2\pi} \int_{-C} \frac{\Delta}{(w - i\Delta)(\gamma w - ia)} \, dw = 0,$$

where the previous integral has again been evaluated using residue calculus. As a result, we obtain that $f(1) = 0$. Applying these considerations with $a = \Delta$ and $a = 2\phi'(\theta_1)$ shows that the integral in (49) equals zero. We also can apply residue calculus to evaluate (48) directly as follows. Consider

$$f_1(w) = \frac{\Delta}{w(w - i\Delta)} \log\left(\frac{2(\gamma(\theta_1) - \gamma(\theta_1 + iw))}{w(w - i\Delta)}\right).$$

Using the change of variables $w = h + i\Delta$ and the definition of $\Delta = \theta_1 - \theta_0$ with $\gamma(\theta_1) = \gamma(\theta_0)$, we can evaluate the residue of $f_1$ at $w = i\Delta$ as Residue$(f_1; i\Delta) = -i\log(-2\gamma'(\theta_0)/\Delta)$. We also can obtain Residue$(f_1; 0) = i\log(2\gamma'(\theta_1)/\Delta)$. Therefore, using residue calculus, we obtain that the integral in (48) equals

$$-i\log(-2\gamma'(\theta_0)/(2\gamma'(\theta_1))) = -i\log(\gamma'(\theta_1)/\gamma'(\theta_1)) = 0,$$

since in the case of symmetric distributions $\gamma'(\lambda)$ is odd and $\theta_1 = -\theta_0$.

Finally, we analyze $B(\Delta)$. Note that

(50) $$B(\Delta) = \frac{1}{2\pi} \int_{|\lambda| \geq \delta} \frac{-\Delta}{(\Delta + i\lambda)i\lambda} \log(2(1 - g(\lambda))\lambda^{-2}) \, d\lambda$$

(51) $$+ \frac{1}{2\pi} \int_{|\lambda| \geq \delta} \frac{-\Delta}{(\Delta + i\lambda)i\lambda} \log\left(1 - \frac{\lambda H(\theta_1, \lambda)}{\lambda - 2i\phi'(\theta_1)}\right) d\lambda.$$

Let $B_1(\Delta)$ and $B_2(\Delta)$ be defined as (50) and (51), respectively. Since $X_1$ is symmetric, it follows that $\log(2(1-g(\lambda))\lambda^{-2})$ is real. As a result, we obtain, just by integrating the real and imaginary parts of the integrand in $B_1(\Delta)$,

$$(52) \qquad B_1(\Delta) = \frac{1}{2\pi} \int_{|\lambda| \geq \delta} \frac{\Delta}{\Delta^2 + \lambda^2} \log(2(1-g(\lambda))\lambda^{-2})\, d\lambda.$$

Expression (52) yields an asymptotic expansion in odd powers of $\Delta$ for $B_1(\Delta)$. Again, integrating the real and imaginary parts in $B_2(\Delta)$, we obtain

$$(53) \qquad B_2(\Delta) = \frac{1}{2\pi} \int_{|\lambda| \geq \delta} \frac{\Delta}{(\Delta^2 + \lambda^2)} \operatorname{Re} \log\left(1 - \frac{\lambda H(\theta_1, \lambda)}{\lambda - 2i\phi'(\theta_1)}\right) d\lambda$$

$$(54) \qquad \qquad - \frac{1}{2\pi} \int_{|\lambda| \geq \delta} \frac{\Delta^2}{(\Delta^2 + \lambda^2)\lambda} \operatorname{Im} \log\left(1 - \frac{\lambda H(\theta_1, \lambda)}{\lambda - 2i\phi'(\theta_1)}\right) d\lambda.$$

The previous identity for $B_2(\Delta)$ is obtained by observing that the integral of the imaginary part must vanish. This occurs because, for all $\theta_1$ small, the function $\log(1 - i\lambda H(\theta_1, \lambda)/(i\lambda + 2\phi'(\theta_1)))$ satisfies the parity property, which can be verified by observing that, since $\gamma(i\lambda) = E_0 \cos(\lambda X) + iE_0 \sin(\lambda X)$, it follows that $h_k(i\lambda) \in \mathcal{P}$; also, using Proposition 6, we obtain that $i\lambda/(i\lambda + 2\phi'(\theta_1))$ satisfies the parity property. Therefore, the closure properties proved in Proposition 6 together with an expansion of the logarithm yield that $\log(1 - i\lambda H(\theta_1, \lambda)/(i\lambda + 2\phi'(\theta_1))) \in \mathcal{P}$, which justifies (53) and (54). For notational convenience, let us define

$$(55) \qquad C(\theta_1, \lambda) = \sum_{k=1}^{\infty} h_{2k}(i\lambda)\theta_1^{2k}/(2k)!$$

and

$$(56) \qquad D(\theta_1, \lambda) = -i \sum_{k=1}^{\infty} h_{2k-1}(i\lambda)\theta_1^{2k-1}/(2k-1)!,$$

where $h_k(i\lambda) = (\gamma^{(k)}(i\lambda) - \mu_k)/(1 - \gamma(i\lambda)) - 2i\mu_{k+1}/\lambda$. Since the distribution of $X_1$ is symmetric, we have that $\gamma(i\lambda)$ is even and real. Moreover, we also have that $h_k(i\lambda)$ is even if and only if $k$ is even. We also can see that $\operatorname{Re}(H(\theta, \lambda)) \triangleq C(\theta, \lambda)$ and $\operatorname{Im}(H(\theta, \lambda)) \triangleq D(\theta, \lambda)$ are even and odd functions of both $\theta$ and $\lambda$ [meaning that, for every $\theta \in (-\eta/2, \eta/2)$ fixed, $C(\theta, \cdot)$ is even and, similarly, for each $\lambda \in \mathbb{R}$, $C(\cdot, \lambda)$ is also even on $(-\eta/2, \eta/2)$, say]. Using this notation, we can write

$$(57) \qquad \frac{\lambda H(\theta_1, \lambda)}{\lambda - 2\phi'(\theta_1)i} = \frac{\lambda^2 C(\theta_1, \lambda) - 2\phi'(\theta_1)\lambda D(\theta_1, \lambda)}{\lambda^2 + (2\phi'(\theta_1))^2}$$

$$(58) \qquad \qquad + i\frac{2\phi'(\theta_1)\lambda C(\theta_1, \lambda) + \lambda^2 D(\theta_1, \lambda)}{\lambda^2 + (2\phi'(\theta_1))^2}.$$

COMPLETE DIFFUSION APPROXIMATIONS 33

Let us define $\overline{C}(\theta_1, \lambda)$ and $\overline{D}(\theta_1, \lambda)$ as the real and imaginary parts of $\lambda H(\theta_1, \lambda)/(\lambda - 2\phi'(\theta_1)i)$, respectively, as indicated in the corresponding expressions (57) and (58). Since $\lambda H(\theta_1, \lambda)/(\lambda - 2\phi'(\theta_1)i)$ holds the parity property, $\overline{C}(\theta_1, \lambda)$ and $\overline{D}(\theta_1, \lambda)$ are even and odd function in both arguments $\theta_1$ and $\lambda$. By symmetry of the distribution of $X_1$, we have that $\Delta = 2\theta_1$, also as a consequence of symmetry, $2\phi'(\theta_1)$ is an odd (real) analytic function of $\theta_1$ at the origin, which implies that $(2\phi'(\theta_1))^2$ is even. Hence, using the expansion of $\log(1 - z)$ at $z = 0$ in expressions (53) and (54) (justified by virtue of Proposition 3), we see that an asymptotic expansion for the integral (53) involves expanding expressions of the form

$$(59) \qquad K(\theta_1) = \frac{1}{2\pi} \int_{|\lambda| \geq \delta} \frac{\Delta}{(\Delta^2 + \lambda^2)} \overline{C}(\theta_1, \lambda)^k \overline{D}(\theta_1, \lambda)^{2m} \, d\lambda,$$

where $K(\theta_1)$ is an even function of $\theta_1$ which is also (real) analytic at the origin. This implies, in view of (55) to (58) and the properties of $2\phi'(\theta_1)$ discussed before, that an asymptotic expansion for (59) must be given in odd powers of $\Delta$ only, which must be also the case for the integral in (53). The treatment for the integral (54) is completely analogous and also yields an asymptotic expansion in odd powers of $\Delta$. This yields the conclusion of the theorem. $\square$

**Acknowledgments.** We thank Professor David Siegmund for insightful conversations on the topic of the present paper and also an Associate Editor and the referees for their valuable comments.## REFERENCES

[1] ASMUSSEN, S. (2001). *Ruin Probabilities*. World Scientific, Singapore. MR1794582
[2] ASMUSSEN, S. (2003). *Applied Probability and Queues*. Springer, New York. MR1978607
[3] BROADIE, M., GLASSERMAN, P. and KOU, S. (1997). A continuity correction for discrete barrier options. *Math. Finance* **7** 325–349. MR1482707
[4] BREIMAN, L. (1992). *Probability*. Addison–Wesley, Reading, MA. MR0229267
[5] BUTZER, P. and NESSEL, R. (1971). *Fourier Analysis and Approximation* **1**. Birkhäuser, Boston. MR0510857
[6] CHANG, J. (1992). On moments of the first ladder height of random walks with small drift. *Ann. Appl. Probab.* **2** 714–738. MR1177906
[7] CHANG, J. and PERES, Y. (1997). Ladder heights, Gaussian random walks and the Riemann zeta function. *Ann. Probab.* **25** 787–802. MR1434126
[8] GLASSERMAN, P. and LIU, T. (1997). Corrected diffusion approximations for multistage production-inventory systems. *Math. Oper. Res.* **12** 186–201. MR1436579
[9] KIEFER, J. and WOLFOWITZ, J. (1956). On the characteristics of the general queueing process, with applications to random walk. *Ann. Math. Statist.* **27** 147–161. MR0077019
[10] KINGMAN, J. (1963). Ergodic properties of continuous time Markov processes and their discrete skeletons. *Proc. London. Math. Soc.* **13** 593–604. MR0154334

[11] LAI, T. (1976). Asymptotic moments of random walks with applications to ladder variables and renewal theory. *Ann. Probab.* **4** 51–66. MR0391265
[12] LINDLEY, D. (1952). The theory of a queue with a single-server. *Proc. Cambridge Philos. Soc.* **48** 277–289. MR0046597
[13] LOTOV, V. (1996). On some boundary crossing problems for Gaussian random walks. *Ann. Probab.* **24** 2154–2171. MR1415246
[14] RUDIN, W. (1987). *Real and Complex Analysis.* McGraw–Hill, New York. MR0924157
[15] SIEGMUND, D. (1979). Corrected diffusion approximations in certain random walk problems. *Adv. in Appl. Probab.* **11** 701–719. MR0544191
[16] SIEGMUND, D. (1985). *Sequential Analysis.* Springer, New York. MR0799155
[17] WOODROOFE, M. (1979). Repeated likelihood ratio tests. *Biometrika* **66** 453–463. MR0556732

DEPARTMENT OF STATISTICS
1 OXFORD ST. SCIENCE CENTER, 7TH FLOOR
HARVARD UNIVERSITY
CAMBRIDGE, MASSACHUSETTS 02138
USA
E-MAIL: blanchet@fas.harvard.edu

MANAGEMENT SCIENCE AND ENGINEERING
TERMAN ENGINEERING CENTER, 3RD FLOOR
STANFORD UNIVERSITY
STANFORD, CALIFORNIA 94305
USA